\documentclass{amsart}

\usepackage{amsmath,amsthm,amsfonts}

\theoremstyle{plain}
\newtheorem{theorem}{Theorem}[section]
\newtheorem{conjecture}{Conjecture}[section]
\newtheorem{lemma}{Lemma}[section]
\newtheorem{corollary}{Corollary}[section]
\numberwithin{equation}{section}

\newcommand{\qbin}[2]{\genfrac{[}{]}{0pt}{}{#1}{#2}}
\newcommand{\qbins}[2]{{\textstyle\genfrac{[}{]}{0pt}{}{#1}{#2}}}
\newcommand{\Z}{\mathbb{Z}}
\newcommand{\B}{\mathcal{B}}
\newcommand{\I}{\mathcal{I}}
\newcommand{\cf}{\text{cf}}
\newcommand{\sign}{\text{sign}}

\begin{document}

\title[The generalized Borwein conjecture]{The generalized Borwein
conjecture. I. The Burge transform}

\author[Ole Warnaar]{S. Ole Warnaar}

\address{Department of Mathematics and Statistics, 
The University of Melbourne, Vic 3010, Australia}
\email{warnaar@ms.unimelb.edu.au}


\thanks{Work supported by the Australian Research Council}

\begin{abstract}
Given an arbitrary ordered pair of coprime integers $(a,b)$ we obtain 
a pair of identities of the Rogers--Ramanujan type.
These identities have the same product side as the (first) Andrews--Gordon
identity for modulus $2ab\pm 1$, but an altogether different sum side,
based on the representation of $(a/b-1)^{\pm 1}$ as a continued fraction.
Our proof, which relies on the Burge transform, first establishes a
binary tree of polynomial identities. Each identity in this Burge tree
settles a special case of Bressoud's generalized Borwein conjecture.
\end{abstract}

\maketitle

\section{Introduction}
Several years ago P.\ Borwein communicated the following
observation to G.\ E.\ Andrews \cite{Andrews95}.
\begin{conjecture}[First Borwein conjecture]\label{Borconj}
The polynomials $A_n(q)$, $B_n(q)$ and $C_n(q)$ defined by
\begin{equation}\label{Bc}
\prod_{k=1}^n(1-q^{3k-2})(1-q^{3k-1})=A_n(q^3)-qB_n(q^3)-q^2C_n(q^3)
\end{equation}
have nonnegative coefficients.
\end{conjecture}
As so often in mathematics, the simplicity of the above claim is rather 
deceptive, and the conjecture still lacks proof.
Following the motto `if you can't prove it, generalize it',
Andrews \cite{Andrews95}, Bressoud \cite{Bressoud96},
and Ismail, Kim and Stanton \cite{IKS99} have extended the Borwein conjecture
in several directions. This is the first of a series of papers
devoted to Bressoud's generalization.
To see how Bressoud's conjecture arises most naturally from 
Conjecture~\ref{Borconj} we follow Andrews \cite{Andrews95} and 
rewrite $A_n, B_n$ and $C_n$ as a sum over $q$-binomial coefficients.
First we need the usual definitions of the $q$-shifted factorial
\begin{equation*}
(a;q)_n=(a)_n=\prod_{k=0}^{n-1} (1-aq^k) \quad\text{for $n\in\Z_+$}
\end{equation*} 
and $q$-binomial coefficient
\begin{equation}\label{qbin}
\qbin{n}{m}_q=\qbin{n}{m}=
\begin{cases}\displaystyle \frac{(q)_{n}}{(q)_m(q)_{n-m}} &
\text{for $m,n-m\in\Z_+$} \\[3mm]
0 & \text{otherwise.} \end{cases}
\end{equation}
Important later will be that for $m$ and $n-m$ nonnegative integers,
the $q$-binomial coefficient is a polynomial in $q$ 
with only positive coefficients.
Also introducing the notation $(a_1,\dots,a_k;q)_n=
(a_1;q)_n\cdots (a_k;q)_n$, we can apply the $q$-binomial theorem
\cite[Eq. (II.4)]{GR90} to expand the left-hand side of \eqref{Bc} as
\begin{align*}
(q,q^2;q^3)_n&=
\sum_{j=-n}^{n}(-1)^j q^{j(3j-1)/2}\qbin{2n}{n-j}_{q^3} \\
&=\sum_{\mu=-1}^1 (-1)^{\mu} q^{\mu(3\mu-1)/2}
\sum_{j\in\Z}(-1)^j q^{3j(9j+6\mu-1)/2} \qbin{2n}{n-3j-\mu}_{q^3}.
\end{align*}
From this one can read off 
\begin{equation}\label{An}
A_n(q)=\sum_{j\in\Z}(-1)^j q^{j(9j-1)/2}\qbin{2n}{n-3j}
\end{equation}
and similar expressions for $B_n$ and $C_n$.
To pass from this to Bressoud's conjecture we need to
recall an important result from partition theory.

Let $\lambda$ be a partition and $\lambda'$ its conjugate.
The $(i,j)$th node of $\lambda$ is the node in the $i$th row and $j$th
column of the Ferrers diagram of $\lambda$. The $d$th diagonal of $\lambda$
is formed by the nodes with coordinates $(i,i-d)$. The hook difference at
node $(i,j)$ is defined as $\lambda_i-\lambda_j'$.
In \cite[Thm. 1]{ABBBFV87} Andrews \textit{et al.} prove the following
theorem using recurrences.
\begin{theorem}\label{thmhp}
The generating function $D_{K,i}(N,M;\alpha,\beta)$ of partitions with at
most $M$ parts, largest part not exceeding $N$, and hook differences on the
$(1-\beta)$th diagonal at least $\beta-i+1$ and on the $(\alpha-1)$th
diagonal at most $K-\alpha-i-1$ is given by
\begin{multline}\label{hookp}
D_{K,i}(N,M;\alpha,\beta;q)=D_{K,i}(N,M;\alpha,\beta) \\[1mm]
=\sum_{j\in\Z} \biggr\{q^{j((\alpha+\beta)Kj+K\beta-(\alpha+\beta)i)}
\qbin{M+N}{M-Kj}-
q^{((\alpha+\beta)j+\beta)(Kj+i)}\qbin{M+N}{M-Kj-i}\biggl\}.
\end{multline}
Here the following conditions apply:
$\alpha,\beta\geq 0$ and $\beta-i\leq N-M\leq K-\alpha-i$
with the added restrictions that the largest part exceeds $M-i$
if $\beta=0$ and the number of parts exceeds $N+i$ if $\alpha=0$.
\end{theorem}
If we follow \cite{Bressoud96} and define
$G(N,M;\alpha,\beta,K)=D_{2K,K}(N,M;\alpha,\beta)$, then
\begin{equation}\label{Gdef}
G(N,M;\alpha,\beta,K)=
\sum_{j\in\Z}(-1)^j q^{\frac{1}{2}Kj((\alpha+\beta)j+\alpha-\beta)}
\qbin{M+N}{N-Kj}.
\end{equation}
Now observe that the expression \eqref{An} for $A_n(q)$ is exactly 
of this type. Explicitly,
\begin{align*}
A_n(q)&=G(n,n;4/3,5/3,3)\\
B_n(q)&=G(n+1,n-1;2/3,7/3,3)\\
C_n(q)&=G(n+1,n-1;1/3,8/3,3).
\end{align*}
Unfortunately there is the complication that $\alpha$ and $\beta$ have
assumed noninteger values so that Theorem~\ref{thmhp} cannot be applied
to interpret $A_n$--$C_n$ as generating functions.
Although no progress has been made in proving Conjecture~\ref{Borconj},
the following generalization is clearly suggested \cite{Bressoud96}.
\begin{conjecture}[Bressoud's generalized Borwein conjecture]\label{conjBr}
Let $K$ be a positive integer and $N,M,\alpha K,\beta K$
be nonnegative integers such that $1\leq \alpha+\beta\leq 2K-1$ (strict
inequalities when $K=2$) and $\beta-K\leq N-M\leq K-\alpha$.
Then $G(N,M;\alpha,\beta,K;q)$ is a polynomial in $q$ with
nonnegative coefficients.
\end{conjecture}
Of course, when both $\alpha$ and $\beta$ are integers the conjecture
becomes a special case of Theorem~\ref{thmhp}. 
When $M+N$ is even and $\alpha=(K-N+M\pm 1)/2$, $\beta=(K+N-M\mp 1)/2$
the conjecture was proven by Ismail \textit{et al.} \cite[Thm. 5]{IKS99}.
We should also remark that not all of the polynomials
$G(N,M;\alpha,\beta,K;q)$ are independent. Besides the obvious symmetry
\begin{equation*}
G(N,M;\alpha,\beta,K)=G(M,N;\beta,\alpha,K)
\end{equation*}
there also holds
\begin{multline}\label{Ginv}
G(N,M;\alpha,\beta,K;1/q) \\ =q^{-MN}G(N,M;K-\alpha-N+M,K-\beta+N-M,K;q)
\end{multline}
and
\begin{align*}
G(N,M;\alpha,\beta,K)
&=G(N,M-1;\alpha,\beta,K)+q^M G(N-1,M;\alpha+1,\beta-1,K) \\
&=G(N-1,M;\alpha,\beta,K)+q^N G(N,M-1;\alpha-1,\beta+1,K)
\end{align*}
as follows from 
\begin{equation}\label{qbininv}
\qbin{n}{m}_{1/q}=q^{m(m-n)}\qbin{n}{m}_q
\end{equation}
and
\begin{equation}\label{qbinrec}
\qbin{n}{m}=\qbin{n-1}{m}+q^{n-m}\qbin{n-1}{m-1}=
\qbin{n-1}{m-1}+q^m \qbin{n-1}{m}.
\end{equation}
These symmetries, not all of which are independent, are consistent
with the bounds imposed on Conjecture~\ref{conjBr}.

This paper is the first in a series in which we apply known 
(the present paper) and new (subsequent papers) transformation
formulas for simple $q$-polynomials
to obtain identities for $G(N,M;\alpha,\beta,K)$ 
that prove its positivity of coefficients.
A few simple examples of such identities can already be found in the 
literature and we quote
\begin{align}
G(L,L;1/2,1,2)&=\sum_{n=0}^L q^{n L}\qbin{L}{n} \label{GB1} \\
G(L,L;1,3/2,2)&=\sum_{n=0}^L q^{n^2}\qbin{L}{n} \label{GB2} \\[2mm]
G(L,L;1/2,3/2,2)&=(1+q^L)(-q^2;q^2)_{L-1} \nonumber .
\end{align}
The first two, which are dual in the sense of \eqref{Ginv},
were found by Bressoud~\cite[Eqs. (8) and (9)]{Bressoud81} and are
bounded analogues of the Euler- and first Rogers--Ramanujan identity.
The third identity is due to Ismail \textit{et al.} \cite[Prop. 2 (3)]{IKS99}.
Although we will add infinitely many new identities to the above three,
we failed to obtain identities that would settle the original Borwein
conjecture. What our results do suggest however is that this is perhaps
more a practical than fundamental problem. 
Often it is necessary to first guess identities before proving
them and it seems that when it comes to the (generalized)
Borwein conjecture the guessing is by far the hardest part of the game.
For example, almost anyone familiar with $q$-series will have
little trouble proving the following identities for
$G(L,L;1,2/3,3)$ and $G(L+1,L-1;1/3,4/3,3)$:
\begin{equation}\label{RR1}
\sum_{j\in\Z}(-1)^j q^{\frac{1}{2}j(5j+1)}\qbins{2L}{L-3j}
=\sum_{n,i\geq 0} q^{n^2+i(L+n)}\qbins{2L-2i-n}{n}\qbins{L-i-n}{i} 
\end{equation}
and
\begin{equation}\label{RR2}
\sum_{j\in\Z}(-1)^j q^{\frac{1}{2}j(5j+3)}\qbins{2L}{L-3j-1}
=\sum_{n,i\geq 0} q^{n(n+1)+i(L+n+1)}\qbins{2L-2i-n-1}{n}\qbins{L-i-n-1}{i},
\end{equation}
but to guess these results is quite a bit harder.
In particular we note that when $L$ tends to infinity
only the terms with $i=0$ contribute to the sums on the right. 
Also using Jacobi's triple product identity \cite[Thm. 2.8]{Andrews76} thus
yields for $|q|<1$
\begin{equation*}
\sum_{n=0}^{\infty}\frac{q^{n^2}}{(q)_n}
=\frac{(q^2,q^3,q^5;q^5)}{(q;q)_{\infty}}
\quad\text{and}\quad
\sum_{n=0}^{\infty}\frac{q^{n(n+1)}}{(q)_n}
=\frac{(q,q^4,q^5;q^5)}{(q;q)_{\infty}},
\end{equation*}
which are the celebrated Rogers--Ramanujan identities \cite{Rogers94}.
Admittedly, we did not guess \eqref{RR1} and \eqref{RR2} but obtained
them by transforming simpler identities, but since we were not so lucky with
the original Borwein conjecture, a good guess is exactly what is needed.
(See also the discussion in section~\ref{secdis}.)

As a bonus of our identities for $G(M,N;\alpha,\beta,K)$ we obtain
many new identities of the Rogers--Ramanujan type.
Some nice examples worth stating in the introduction are the following
four series of identities featuring the Fibonacci numbers $F_k$ defined
recursively as $F_0=0$, $F_1=1$ and $F_k=F_{k-1}+F_{k-2}$.
\begin{theorem}\label{Fib1}
For $|q|<1$ and $k\geq 4$ there holds
\begin{multline*}
\sum_{m_1,\dots,m_{k-2}\geq 0}
\frac{q^{m_1^2+\cdots+m_{k-2}^2}}{(q)_{2m_1}}
\Bigl(\prod_{j=2}^{k-3}\qbins{m_{j-1}+m_j-m_{j+1}}{2m_j}\Bigr)
\qbins{m_{k-3}}{m_{k-2}} \\
=\frac{(q^{F_kF_{k-1}},q^{F_kF_{k-1}+1},q^{2F_kF_{k-1}+1};
q^{2F_kF_{k-1}+1})_{\infty}}{(q;q)_{\infty}}.
\end{multline*}
\end{theorem}
Taking $\sum_{m\geq 0}q^{m^2}/(q)_m$ as the left-hand side when $k=3$
includes the first Rogers--Ramanujan identity in this series of identities.
\begin{theorem}\label{Fib2}
For $|q|<1$ and $k\geq 4$ there holds
\begin{multline*}
\sum_{m_1,\dots,m_{k-2}\geq 0}
\frac{q^{m_1^2+\cdots+m_{k-3}^2+m_{k-3}m_{k-2}}}{(q)_{2m_1}}
\Bigl(\prod_{j=2}^{k-3}\qbins{m_{j-1}+m_j-m_{j+1}}{2m_j}\Bigr)
\qbins{m_{k-3}}{m_{k-2}} \\
=\frac{(q^{F_kF_{k-1}-1},q^{F_kF_{k-1}},q^{2F_kF_{k-1}-1};
q^{2F_kF_{k-1}-1})_{\infty}}{(q;q)_{\infty}}.
\end{multline*}
\end{theorem}
\begin{theorem}\label{Fib3}
For $|q|<1$ and $k\geq 5$ there holds
\begin{multline*}
\sum_{m_1,\dots,m_{k-2}\geq 0}
\frac{q^{(m_1+m_2)^2+m_2^2+\cdots+m_{k-2}^2}}
{(q)_{m_1}(q)_{2m_2}}
\Bigl(\prod_{j=3}^{k-3}\qbins{m_{j-1}+m_j-m_{j+1}}{2m_j}\Bigr)
\qbins{m_{k-3}}{m_{k-2}} \\
=\frac{(q^{F_kF_{k-2}},q^{F_kF_{k-2}+1},q^{2F_kF_{k-2}+1};
q^{2F_kF_{k-2}+1})_{\infty}}{(q;q)_{\infty}}.
\end{multline*}
\end{theorem}
The cases $k=3$ and $k=4$ may be included by taking
$\sum_{m\geq 0}q^{m^2}/(q)_m$ and
$\sum_{m_1,m_2\geq 0}q^{(m_1+m_2)^2+m_2^2}/(q)_{m_1}(q)_{m_2}$ as
respective left-hand sides.
\begin{theorem}\label{Fib4}
For $|q|<1$ and $k\geq 5$ there holds
\begin{multline*}
\sum_{m_1,\dots,m_{k-2}\geq 0}
\frac{q^{(m_1+m_2)^2+m_2^2+\cdots+m_{k-3}^2+m_{k-3}m_{k-2}}}
{(q)_{m_1}(q)_{2m_2}}
\Bigl(\prod_{j=3}^{k-3}\qbins{m_{j-1}+m_j-m_{j+1}}{2m_j}\Bigr)
\qbins{m_{k-3}}{m_{k-2}} \\
=\frac{(q^{F_kF_{k-2}-1},q^{F_kF_{k-2}},q^{2F_kF_{k-2}-1};
q^{2F_kF_{k-2}-1})_{\infty}}{(q;q)_{\infty}}.
\end{multline*}
\end{theorem}

\section{The Burge transform}
The Burge transform, which is a generalization of a special case of the
Bailey lemma, provides an iterative method for proving
polynomial analogues of $q$-series identities \cite{Burge93}.
Before we state the actual transform let us define the polynomial
$\B(L,M,a,b;q)$ as
\begin{equation}\label{Bdef}
\B(L,M,a,b;q)=\B(L,M,a,b)=\qbin{L+M+a-b}{L+a}\qbin{L+M-a+b}{L-a}.
\end{equation}
Note that $\B(L,M,a,b)$ is nonzero iff $L+a$ and $M+b$ are integers
such that $0\leq |b|\leq M$ and $0\leq |a|\leq L$. Moreover, it satisfies
the symmetries
\begin{align}
\B(L,M,-a,-b)&=\B(L,M,a,b)  \notag \\[2mm]
\B(L,M,a,b)&=\B(M,L,b,a)  \label{Bsymm} \\[2mm]
\B(L,M,a,b;1/q)& =q^{2ab-2LM}\B(L,M,a,b;q) \label{Brecip}
\end{align}
and becomes proportional to the $q$-binomial coefficient when
either $L$ or $M$ tends to infinity
\begin{equation}\label{Blim}
\lim_{M\to\infty}\B(L,M,a,b)=\qbin{2L}{L-a}/(q)_{2L}.
\end{equation}
The most interesting properties of $\B$ are however the following two
transformations.
\begin{theorem}\label{thmBurge}
For $L,M,a,b$ integers such that not $-L+a\leq -b\leq L+a<b\leq M$
or $-L-a\leq b\leq L-a<-b\leq M$,
\begin{equation}\label{B1}
\sum_{i\geq 0}q^{i^2}\qbin{2L+M-i}{2L}\B(L-i,i,a,b)=q^{b^2}\B(L,M,a+b,b)
\end{equation}
and
\begin{equation}\label{B2}
\sum_{i\geq 0}q^{i^2}\qbin{2L+M-i}{2L}\B(i,L-i,b,a)=q^{b^2}\B(L,M,a+b,b).
\end{equation}
\end{theorem}
These two results are known as the Burge transform
\cite{Burge93,FLW97,SW00}, and are an immediate consequence of
the $q$-Saalsch\"utz sum \cite[Eq. (II.12)]{GR90}, or, equivalently, 
of \cite[Eq. (3.3.11)]{Andrews76}.
The second transform of course follows from the first by the symmetry 
relation \eqref{Bsymm}.
The conditions imposed on the parameters are due to the fact that the
left side may be zero when the right side is not. They perhaps
appear cumbersome but are in fact quite innocent.
In all our applications of the Burge transform $a$ and $b$ will
have the same signature. It is not hard to see that the
conditions then become void. For example, assuming $0<b\leq M$
we need to inspect the first condition only. But $a\geq 0$
clearly avoids the occurrence of $-b\leq L+a<b\leq L-a$.
A similar reasoning applies for negative $b$.

\section{The generalized Borwein conjecture}\label{secGB}
To derive manifestly positive representations for $G(N,M;\alpha,\beta,K)$
using the Burge transform the following lemma will be crucial.
\begin{lemma}\label{lemB}
For $L,M\geq 0$ there holds
\begin{equation}\label{Bnew}
\sum_{j\in\Z}(-1)^j q^{\frac{1}{2}j(3j+1)}\B(L,M,j,j)=\qbin{L+M}{M}.
\end{equation}
\end{lemma}
\begin{proof}
Take Rogers' $q$-Dougall sum \cite[(II.21)]{GR90}
\begin{equation*}
\sum_{k=0}^n \frac{1-aq^{2k}}{1-a}\frac{(a,b,c,q^{-n};q)_k}
{(q,aq/b,aq/c,aq^{n+1};q)_k}\Bigr(\frac{aq^{n+1}}{bc}\Bigl)^k=
\frac{(aq,aq/bc;q)_n}{(aq/b,aq/c;q)_n}
\end{equation*}
and let $a\to 1$, $b \to q^{-M}$, $c\to\infty$ and $n\to L$.
After some simplifications this gives \eqref{Bnew}.
\end{proof}
When either $L$ or $M$ tends to infinity \eqref{Bnew} simplifies to Rogers'
polynomial analogue of the Euler identity \cite[\S 1]{Rogers17} 
\begin{equation*}
\sum_{j\in\Z}(-1)^j q^{\frac{1}{2}j(3j+1)}\qbin{2M}{M-j}=
\frac{(q)_{2M}}{(q)_M}.
\end{equation*}

In the following we will iterate the two Burge transformations \eqref{B1}
and \eqref{B2} to transform \eqref{Bnew} into a binary tree of polynomial
identities. First, by application of either \eqref{B1} or \eqref{B2} one
finds
\begin{equation}\label{Bnewp}
\sum_{j\in\Z}(-1)^j q^{\frac{1}{2}j(5j+1)}\B(L,M,2j,j)
=\sum_{n\geq 0}q^{n^2}\qbin{2L+M-n}{2L}\qbin{L}{n}.
\end{equation}
This is a doubly-bounded analogue of the first Rogers--Ramanujan
identity. When $M$ goes to infinity we recover Bressoud's identity
\eqref{GB2} and when $L$ goes to infinity we
find the following well-known specialization of Watson's $_{8}\phi_7$ 
transform (see e.g. \cite[Eq. (1.11)]{Bressoud81} or \cite[Eq. (39)]{Paule85})
\begin{equation*}
\sum_{j\in\Z}(-1)^j q^{\frac{1}{2}j(5j+1)}\qbin{2M}{M-j}
=\frac{(q)_{2M}}{(q)_M}\sum_{n\geq 0}q^{n^2}\qbin{M}{n}.
\end{equation*}

Without too much effort one can utilize \eqref{Bnewp} to also obtain
a doubly-bounded analogue of the second Rogers--Ramanujan identity.
\begin{lemma}
For $L,M\geq 0$ there holds
\begin{equation}\label{comp}
\sum_{j\in\Z}(-1)^j q^{\frac{1}{2}j(5j+3)}
\qbins{L+M+j}{M-j-1}\qbins{L+M-j}{M+j}
=\sum_{n\geq 0}q^{n(n+1)}\qbins{2L+M-n}{2L+1}\qbins{L}{n}.
\end{equation}
\end{lemma}
\begin{proof}
Suppressing their $L$-dependence we denote the left sides of \eqref{comp} 
and \eqref{Bnewp} by $f_M$ and $g_M$, respectively.
By application of the recurrence \eqref{qbinrec}
\begin{multline*}
f_M=\sum_{j\in\Z}(-1)^j q^{\frac{1}{2}j(5j+3)}
\qbins{L+M+j}{M-j-1}\qbins{L+M-j-1}{M+j-1} \\
+q^M\sum_{j\in\Z}(-1)^j q^{\frac{5}{2}j(j+1)}
\qbins{L+M+j}{M-j-1}\qbins{L+M-j-1}{M+j}.
\end{multline*}
Since the second term on the right changes sign after the variable change
$j\to -j-1$ it vanishes. The first term is again split using
\eqref{qbinrec} leading to 
\begin{equation}\label{frec}
f_M=f_{M-1}+q^{M-1}g_{M-1}.
\end{equation}
Next we let $f_M$ and $g_M$ denote the right sides of \eqref{comp}
and \eqref{Bnewp}. A single application of \eqref{qbinrec} show that
\eqref{frec} again holds.
Since equation \eqref{comp} trivializes to $0=0$ for $M=0$
this settles the lemma.
\end{proof}

It is easy to see that by computing $\eqref{comp}+q^M \eqref{Bnewp}$
we also get
\begin{equation}\label{comp2}
\sum_{j\in\Z}(-1)^j q^{\frac{1}{2}j(5j+3)}\qbins{L+M+j+1}{M-j}
\qbins{L+M-j}{M+j}
=\sum_{n\geq 0}q^{n(n+1)}\qbins{2L+M-n+1}{2L+1}\qbins{L}{n}.
\end{equation}
Although the identities \eqref{comp} and \eqref{comp2} cannot be written
in terms of the polynomials $\B(L,M,a,b)$ they can be iterated by
a simple modification of the Burge transform. We will not pursue
this here and only take the large $M$ limit to arrive at the isolated
positivity result
\begin{equation*}
G(L+1,L;1/2,2,2)=
\sum_{j\in\Z}(-1)^j q^{\frac{1}{2}j(5j+3)}\qbin{2L+1}{L-2j}
=\sum_{n\geq 0}q^{n(n+1)}\qbin{L}{n}.
\end{equation*}

After this intermezzo we continue to iterate \eqref{Bnewp}.
To state the general result some more notation is needed.
Assume that $(a,b)$ is a pair of coprime integers such that $1\leq b<a$,
and define a nonnegative integer $n$ and positive integers $a_0,\dots,a_n$
as the order and partial quotients of the continued fraction representation
of $(a/b-1)^{\sign(a-2b)}$ ($\sign(0)=0$), i.e.,
\begin{equation}\label{cfe}
\Bigl(\frac{a}{b}-1\Bigr)^{\sign(a-2b)}=[a_0,\dots,a_n]=
a_0+\frac{1}{a_1+\cfrac{1}{\cdots+\cfrac{1}{a_n}}}.
\end{equation}
We denote the continued fraction corresponding to $(a,b)$ by $\cf(a,b)$,
and note the obvious symmetry $\cf(a,b)=\cf(a,a-b)$.
By abuse of notation we sometimes write $\cf(a,b)=(a/b-1)^{\sign(a-2b)}$.

Before we continue, we make the following important remark concerning continued
fractions. For any admissible $(a,b)$ such that 
$(a,b)\neq (2,1)$ the continued fraction $\cf(a,b)$ is not unique. Indeed,
for $c_n\geq 2$ the continued fractions $[c_0,\dots,c_n]$ and 
$[c_0,\dots,c_n-1,1]$ represent the same rational number.
This means in particular that the order and partial quotients 
of $\cf(a,b)$ defined in \eqref{cfe} are not unique. 
In the following we will use these to define several other quantities, 
which are therefore not unique either.
Later we will show however, that whichever choice is made, the final
objects of interest are independent of the choice of representation.
This allows use to choose at our convenience the representation with either
$a_n=1$ or $a_n\geq 2$. For $(a,b)=(2,1)$ we of course
only have $\cf(2,1)=[1]$.

We now continue by further defining the partial sums 
$t_j=\sum_{k=0}^{j-1}a_k$ for $j=1,\dots,n+1$. 
We also introduce $t_0=0$ and $d(a,b)=t_{n+1}$. Note that $d(a,b)$ is
insensitive to the choice of representation of $\cf(a,b)$.
Finally we define a $d(a,b)\times d(a,b)$ matrix $\I(a,b)$ with entries
\begin{equation*}
\I(a,b)_{j,k}=\begin{cases}
\delta_{j,k+1}+\delta_{j,k-1} & \text{for $j\neq t_i$} \\
\delta_{j,k+1}+\delta_{j,k}-\delta_{j,k-1} & \text{for $j=t_i$}
\end{cases}
\end{equation*}
and a corresponding Cartan-type matrix $C(a,b)=2I-\I(a,b)$ where $I$ is
the $d(a,b)\times d(a,b)$ identity matrix.
Note that the matrix $\I(a,b)$ has the following block-structure:
\begin{equation*}
\I(a,b)=\left(\begin{array}{ccc|ccc|ccc}
&&&&&&&& \\
&T_{a_0}&&&&&&& \\
&&&-1&&&&& \\ \hline
&&1&&&&&& \\
&&&&\ddots&&&& \\
&&&&&&-1&& \\ \hline
&&&&&1&&& \\
&&&&&&&T_{a_n}& \\
&&&&&&&&
\end{array}\right)
\end{equation*}
where $T_i$ is the incidence matrix of the tadpole graph with $i$ vertices,
i.e., $(T_i)_{j,k}=\delta_{j,k+1} +\delta_{j,k-1}+\delta_{j,k}\delta_{j,i}$.
A change from the representation of $\cf(a,b)$ with $a_n\geq 2$
to the representation with $a_n=1$ corresponds to the transformation
\begin{equation}\label{Ttrafo}
T_{a_n}\to\left(\begin{array}{rrr|r}&&& \\
& T_{a_n-1} && \\ &&& -1 \\ \hline &&1&1 \end{array}\right)
\end{equation}
which leaves all but one matrix elements unchanged.
(Given that $\I(a,b)$ and $C(a,b)$ depend not only on $(a,b)$ but
also on the choice of representation of $\cf(a,b)$ the fastidious reader 
might prefer the notation $\I(\cf(a,b))$ and $C(\cf(a,b))$.)

With the above notation we define a polynomial for each ordered pair of
positive, coprime integers $(a,b)$ by
\begin{equation}\label{Fab}
F_{a,b}(L,M)=\sum_{m\in\Z_+^{d(a,b)}}q^{mC(a,b)m}\qbins{2L+M-m_1}{2L}
\prod_{j=1}^{d(a,b)}\qbins{\tau_j m_j+n_j}{\tau_j m_j}
\end{equation}
for $a\leq 2b$, and
\begin{equation}\label{Faab}
F_{a,b}(L,M)=\sum_{m\in \Z_+^{d(a,b)}}q^{L(L-2m_1)+mC(a,b)m}
\qbins{L+M+m_1}{2L}
\prod_{j=1}^{d(a,b)}\qbins{\tau_j m_j+n_j}{\tau_j m_j}
\end{equation}
for $a\geq 2b$. Here 
\begin{equation}\label{mCm}
mC(a,b)m=\sum_{j,k=1}^{d(a,b)}m_j C(a,b)_{j,k}m_k=\sum_{j=0}^n
\Bigl(m_{t_j+1}^2+\sum_{k=t_j+1}^{t_{j+1}-1}(m_{k}-m_{k+1})^2\Bigr)
\end{equation}
and $\tau_1=\dots=\tau_{d(a,b)-1}=2$, $\tau_{d(a,b)}=1$.
The auxiliary variables $n_j$ in the summand are integers
defined by the $(m,n)$-system
\begin{equation}\label{mn}
n_j=L\delta_{j,1}-\sum_{k=1}^{d(a,b)}C(a,b)_{j,k} m_k
\quad \text{for $j=1,\dots,d(a,b)$.}
\end{equation}
We remark that the polynomials $F_{a,b}(L,M)$ bear close resemblance
to the much-studied fermionic polynomials arising in
the quasiparticle description of $c<1$ conformal field theory
\cite{BM96,BMS98,FLW97,SW98}.

For each admissible $(a,b)$ the polynomial $F_{a,b}(L,M)$ is defined
twice. To show that these definitions are consistent we
first assume that $(a,b)\neq (2,1)$ and that we have chosen $\cf(a,b)$ such
that $a_n\geq 2$. Making the variable change 
$m_{d(a,b)}\leftrightarrow n_{d(a,b)}$ then yields the representation
with $a_n=1$.
The easiest way to see this is perhaps to first eliminate the $n_j$, then to
make the variable change $m_{d(a,b)}\to m_{d(a,b)-1}-m_{d(a,b)}$ and to 
rewrite this again in the form using the $n_j$. Recalling \eqref{Ttrafo}
and the definition of $T_i$ and $C(a,b)$ this gives the desired result.
When $(a,b)=(2,1)$ there is only one continued fraction representation,
but $F_{2,1}(L,M)$ is defined in both \eqref{Fab} and \eqref{Faab}.
Since $\cf(2,1)=[1]$, one finds $d(2,1)=1$ and $\I(2,1)=C(2,1)=(1)$.
From \eqref{mn} we further find $n_1=L-m_1$.
Hence, according to \eqref{Fab}
\begin{equation*}
F_{2,1}(L,M)=\sum_{m_1\geq 0} q^{m_1^2}\qbin{2L+M-m_1}{2L}\qbin{L}{m_1}
\end{equation*}
and according to \eqref{Faab}
\begin{equation*}
F_{2,1}(L,M)=\sum_{m_1\geq 0} q^{(L-m_1)^2}\qbin{L+M+m_1}{2L}\qbin{L}{m_1}.
\end{equation*}
The variable change $m_1\to L-m_1$ shows these two expressions
are consistent. 

We are now ready to state the identities obtained by applying the
Burge transform to \eqref{Bnewp}.
\begin{theorem}\label{thmmain}
For $L,M$ nonnegative integers and $(a,b)$ a pair of coprime integers
such that $1\leq b<a$ there holds
\begin{equation}\label{eqmain}
\sum_{j\in\Z}(-1)^j q^{\frac{1}{2}j((2ab+1)j+1)}\B(L,M,aj,bj)=F_{a,b}(L,M).
\end{equation}
\end{theorem}
We postpone the proof till the next section and instead
continue with examples and corollaries.

In our first example we take $(a,b)=(7,5)$ and explicitly
calculate $F_{7,5}(L,M)$ and $F_{7,2}(L,M)$. 
As continued fraction we choose the representation
$\cf(7,5)=[2,1,1]$. Hence $n=2$, $a_0=2$, $a_1=a_2=1$, $t_1=2$, $t_2=3$
and $t_3=d(7,5)=4$. The matrices $\I(7,5)$ and $C(7,5)$ are thus found
to be \begin{equation*}
\I(7,5)=\left(\begin{array}{rr|r|r}0&1&0&0 \\ 1&1&-1&0 \\ \hline 
0&1&1&-1 \\ \hline 0&0&1&1 \end{array}
\right), \qquad
C(7,5)=\left(\begin{array}{rr|r|r}2&-1&0&0 \\ -1&1&1&0 \\ \hline
0&-1&1&1 \\ \hline 0&0&-1&1 \end{array}\right).
\end{equation*}
Inserting this in \eqref{Fab} and \eqref{Faab} (using the symmetry
$\cf(7,5)=\cf(7,2)$) yields
\begin{multline*}
F_{7,5}(L,M)=\sum_{m_1,\dots,m_4\geq 0}
q^{m_1^2+(m_1-m_2)^2+m_3^2+m_4^2}
\qbins{2L+M-m_1}{2L}\qbins{L+m_2}{2m_1} \\
\times \qbins{m_1+m_2-m_3}{2m_2}\qbins{m_2+m_3-m_4}{2m_3}\qbins{m_3}{m_4}
\end{multline*}
and
\begin{multline*}
F_{7,2}(L,M)=\sum_{m_1,\dots,m_4\geq 0}
q^{(L-m_1)^2+(m_1-m_2)^2+m_3^2+m_4^2}
\qbins{L+M+m_1}{2L}\qbins{L+m_2}{2m_1} \\
\times \qbins{m_1+m_2-m_3}{2m_2}\qbins{m_2+m_3-m_4}{2m_3}\qbins{m_3}{m_4}.
\end{multline*}

Next we consider Theorem~\ref{thmmain} in the large $M$ limit,
and with the same notation as above we define $F_{a,b}(L)=(q)_{2L}
\lim_{M\to\infty}F_{a,b}(L,M)$. Explicitly this yields
\begin{equation}\label{FabL}
F_{a,b}(L)=\sum_{m\in\Z_+^{d(a,b)}}
q^{mC(a,b)m}\prod_{j=1}^{d(a,b)}\qbins{\tau_j m_j+n_j}{\tau_j m_j}
\end{equation}
for $a\leq 2b$, and
\begin{equation}\label{FaabL}
F_{a,b}(L)=\sum_{m\in\Z_+^{d(a,b)}}
q^{L(L-2m_1)+mC(a,b)m}\prod_{j=1}^{d(a,b)}\qbins{\tau_j m_j+n_j}{\tau_j m_j}
\end{equation}
for $a\geq 2b$. 
Then, using \eqref{Blim}, the following result arises after letting $M$ 
tend to infinity in \eqref{eqmain}.
\begin{corollary}\label{Corab}
For $L$ a nonnegative integer and $(a,b)$ a pair of coprime integers
such that $1\leq b<a$ there holds 
\begin{equation}\label{GF}
G(L,L;b,b+1/a,a)=
\sum_{j\in\Z}(-1)^j q^{\frac{1}{2}j((2ab+1)j+1)}\qbin{2L}{L-aj}=F_{a,b}(L).
\end{equation}
\end{corollary}
Obviously, $F_{a,b}(L)$ is a polynomial with only nonnegative coefficients.
This leads to our next corollary.
\begin{corollary}\label{corGenBor}
For $(a,b)$ a pair of coprime integers such that $1\leq b<a$ the
polynomial $G(L,L;b,b+1/a,a)$ has nonnegative coefficients.
\end{corollary}

Taking the large $L$ limit of Theorem~\ref{thmmain} is more intricate
due to the $L$-dependent term in the exponent of $q$ in \eqref{Faab}.
To overcome this complication one first has to make a change of variables
expressing the kernel of $F_{a,b}(L,M)$ (with $a\geq 2b$) 
in terms of the variables $n_1,\dots,n_{a_0},m_{a_0+1},\dots,m_{d(a,b)}$
instead of $m_1,\dots,m_{d(a,b)}$.
To achieve this, first note that \eqref{mn} implies
\begin{equation}\label{mn2}
m_j=L-jm_{a_0+1}-\sum_{k=1}^{a_0} \min(j,k) n_k
\quad \text{for $j=1,\dots,a_0$,}
\end{equation}
and hence $m_j-m_{j+1}=m_{a_0+1}+n_{j+1}+\cdots+n_{a_0}$.
(If for $(a,b)=(a,1)$ we take the representation $\cf(a,1)=[a-1]$ then
$a_0=d(a,1)=a-1$ in which case $m_{a_0+1}:=0$.) 
Inserting this in \eqref{Faab}, using \eqref{mCm}
and defining $N_j=n_j+\cdots+n_{a_0}$, gives
\begin{multline*}
F_{a,b}(L,M)=\sum_{\substack{n_1,\dots,n_{a_0}\geq 0 \\ 
m_{a_0+1},\dots,m_{d(a,b)}\geq 0}}
q^{(N_1+m_{a_0+1})^2+\cdots+(N_{a_0}+m_{a_0+1})^2} \\
\times q^{\sum_{j,k=a_0+1}^{d(a,b)} m_j C(a,b)_{j,k} m_k}
\qbins{L+M+m_1}{2L}\prod_{j=1}^{a_0}\qbins{\tau_j m_j+n_j}{n_j}
\prod_{j=a_0+1}^{d(a,b)}\qbins{\tau_j m_j+n_j}{\tau_j m_j},
\end{multline*}
with $a\geq 2b$. Here it is to be understood that the auxiliary variables
$m_1,\dots,m_{a_0}$ and $n_{a_0+2},\dots,n_{d(a,b)}$
follow from \eqref{mn2} and \eqref{mn}, respectively.
The auxiliary variable $n_{a_0+1}$ is special in that it needs
both equations for its computation, 
\begin{equation*}
n_{a_0+1}=L-a_0 m_{a_0+1}-\sum_{k=1}^{a_0} k n_k
-\sum_{k=a_0+1}^{d(a,b)}C(a,b)_{a_0+1,k}m_k.
\end{equation*}

After these preliminaries we are prepared to take the large $L$ limit of
Theorem~\ref{thmmain}.
All that remains to be done is to define the polynomials 
$\tilde{F}_{a,b}(M)=(q)_{2M}\lim_{L\to\infty}F_{a,b}(L,M)$, i.e.,
\begin{equation}\label{FabM}
\tilde{F}_{a,b}(M)=
\sum_{m\in\Z^{d(a,b)}_+}
\frac{(q)_{2M}\:q^{mC(a,b)m}}{(q)_{M-m_1}(q)_{\tau_1 m_1}}
\prod_{j=2}^{d(a,b)}\qbins{\tau_j m_j+n_j}{\tau_j m_j}
\end{equation}
for $a\leq 2b$, and
\begin{multline}\label{FaabM}
\tilde{F}_{a,b}(M)=
\sum_{\substack{n_1,\dots,n_{a_0}\geq 0 \\
m_{a_0+1},\dots,m_{d(a,b)}\geq 0}}
\frac{(q)_{2M}\:q^{(N_1+m_{a_0+1})^2+\cdots+(N_{a_0}+m_{a_0+1})^2}}
{(q)_{M-N_1-m_{a_0+1}}
(q)_{n_1}\cdots(q)_{n_{a_0}}(q)_{\tau_{a_0+1}m_{a_0+1}}} \\
\times q^{\sum_{j,k=a_0+1}^{d(a,b)} m_j C(a,b)_{j,k} m_k}
\prod_{j=a_0+2}^{d(a,b)}\qbins{\tau_j m_j+n_j}{\tau_j m_j}
\end{multline}
for $a\geq 2b$. Of course $\tau_{a_0+1}$ is always $2$ except when
$(a,b)=(a,1)$ with the choice $\cf(a,1)=[a-2,1]$. 

Our next corollary can at last be stated as follows.
\begin{corollary}\label{corM}
For $M$ a nonnegative integer and $(a,b)$ a pair of coprime integers
such that $1\leq b<a$ there holds
\begin{equation*}
G(M,M;a,a+1/b,b)=\sum_{j\in\Z}(-1)^j q^{\frac{1}{2}j((2ab+1)j+1)}
\qbin{2M}{M-bj}=\tilde{F}_{a,b}(M).
\end{equation*}
\end{corollary}

Continuing our previous example for $(a,b)=(7,5)$ we find
\begin{equation*}
\tilde{F}_{7,5}(M)=(q)_{2M}\sum_{m_1,\dots,m_4\geq 0}
\frac{q^{m_1^2+(m_1-m_2)^2+m_3^2+m_4^2}}{(q)_{M-m_1}(q)_{2m_1}}
\qbins{m_1+m_2-m_3}{2m_2}\qbins{m_2+m_3-m_4}{2m_3}\qbins{m_3}{m_4}
\end{equation*}
\and
\begin{equation*}
\tilde{F}_{7,2}(M)=(q)_{2M}\sum_{n_1,n_2,m_3,m_4\geq 0}
\frac{q^{(n_1+n_2+m_3)^2+(n_2+m_3)^2+m_3^2+m_4^2}}
{(q)_{M-m_3-n1-n_2}(q)_{n_1}(q)_{n_2}(q)_{2m_3}}\qbins{m_3}{m_4}.
\end{equation*}

Next we return to Corollary~\ref{corM} and take $M$ to infinity.
By the Jacobi triple product identity \cite[Thm. 2.8]{Andrews76}
we then obtain the following two results.
\begin{theorem}\label{RRab}
For $|q|<1$ and $(a,b)$ a pair of coprime integers such that $1\leq b<a$
there holds
\begin{equation*}
\sum_{m\in\Z^{d(a,b)}_+}\frac{q^{mC(a,b)m}}{(q)_{\tau_1 m_1}}
\prod_{j=2}^{d(a,b)}\qbins{\tau_j m_j+n_j}{\tau_j m_j}
=\frac{(q^{ab},q^{ab+1},q^{2ab+1};q^{2ab+1})_{\infty}}{(q;q)_{\infty}}
\end{equation*}
for $a\leq 2b$, and
\begin{multline*}
\sum_{\substack{n_1,\dots,n_{a_0}\geq 0 \\
m_{a_0+1},\dots,m_{d(a,b)}\geq 0}}
\frac{q^{(N_1+m_{a_0+1})^2+\cdots+(N_{a_0}+m_{a_0+1})^2}}
{(q)_{n_1}\cdots(q)_{n_{a_0}}(q)_{\tau_{a_0+1}m_{a_0+1}}} \\
\times q^{\sum_{j,k=a_0+1}^{d(a,b)} m_j C(a,b)_{j,k} m_k}
\prod_{j=a_0+2}^{d(a,b)}\qbins{\tau_j m_j+n_j}{\tau_j m_j}
=\frac{(q^{ab},q^{ab+1},q^{2ab+1};q^{2ab+1})_{\infty}}{(q;q)_{\infty}}
\end{multline*}
for $a\geq 2b$.
\end{theorem}

The simplest case of the theorem is $(a,b)=(k,1)$
or, by symmetry, $(a,b)=(k,k-1)$. Choosing $\cf(k,1)=[k-1]$,
one finds $n=0$, $a_0=k-1$, $d(k,1)=k-1$, and $\I(k,1)=T_{k-1}$.
Hence 
\begin{equation}\label{AGid}
\sum_{n_1,\dots,n_{k-1}\geq 0}
\frac{q^{N_1^2+\cdots+N_{k-1}^2}}{(q)_{n_1}\cdots(q)_{n_{k-1}}}
=\frac{(q^k,q^{k+1},q^{2k+1};q^{2k+1})_{\infty}}{(q;q)_{\infty}}
\end{equation}
for $k\geq 2$ where $N_j=n_j+\cdots+n_{k-1}$, and
\begin{multline*}
\sum_{m_1,\dots,m_{k-1}\geq 0}
\frac{q^{M_1^2+\cdots+M_{k-1}^2}}{(q)_{2m_1}}
\Bigl(\prod_{j=2}^{k-2}\qbins{m_{j-1}+m_{j+1}}{2m_j}\Bigr)
\qbins{m_{k-2}}{m_{k-1}} \\
=\frac{(q^{k(k-1)},q^{k(k-1)+1},q^{2k(k-1)+1};q^{2k(k-1)+1})_{\infty}}
{(q;q)_{\infty}}
\end{multline*}
for $k\geq 3$ with $M_j=m_j-m_{j-1}$ ($m_0=0$).
The first of these two results is of course the (first) Andrews--Gordon
identity for modulus $2k+1$ \cite{Andrews74} which we have now embedded 
in a much larger family of Rogers--Ramanujan-type identities.

Another simple case of Theorem~\ref{RRab} occurs when $(a,b)=(F_k,F_{k-1})$
or $(F_k,F_{k-2})$ where $F_k=((1+\sqrt{5})^k-(1-\sqrt{5})^k)/
(2^k \sqrt{5})$ is the $k$th Fibonacci number;
$F_0=0$, $F_1=1$ and $F_k=F_{k-1}+F_{k-2}$.
Using the recurrence for the Fibonacci numbers and $F_2/F_1=1$ this yields
\begin{equation*}
\cf(F_k,F_{k-1})=\Bigr(\frac{F_k}{F_{k-1}}-1\Bigl)^{-1}=
\frac{F_{k-1}}{F_{k-2}}=1+\frac{1}{\frac{F_{k-2}}{F_{k-3}}}=
\dots =[\underbrace{1,\dots,1}_{k-2}]
\end{equation*}
for $k\geq 3$. Hence $n=k-3$, $a_0=\dots=a_{k-3}=1$ and $d(F_k,F_{k-1})=k-2$.
Inserting this in Theorem~\ref{RRab} and renaming $n_1$ as $m_1$ in the second
equation yields the Theorems~\ref{Fib1} and \ref{Fib3}.

Before we continue to explore the consequences of Theorem~\ref{thmmain}
we should perhaps remark that in the Andrews--Gordon identity \eqref{AGid}
we can freely choose $k$. This means that we can always tune its right-hand
side to coincide with the right-hand side of the first or second
identity of Theorem~\ref{RRab}.
So what we have actually obtained are new representations for the
sum side of the Andrews--Gordon identities. Specifically,
if we consider the Andrews--Gordon identity for modulus $2k+1$ then we
have found new sum sides for each ordered pair of coprime integers $(a,b)$
such that $ab=k$.
The most spectacular aspect of this is perhaps the fact that these
new sums sides, though being considerably more cumbersome to write down,
are much more efficient.
In particular, the Andrews--Gordon identity for modulus $2k+1$ has a 
$(k-1)$-fold sum on the left-hand side whereas for $(a,b)$ such 
that $ab=k$ (and $b>1$)
we obtain a sum side consisting of a $d(a,b)$-fold sum, where
$d(a,b)$ is the sum of the partial quotients in the
continued fraction of $(a/b-1)^{\sign(a-2b)}$.
For instance, from our on-going example we find
\begin{multline*}
\sum_{m_1,\dots,m_4\geq 0}
\frac{q^{m_1^2+(m_1-m_2)^2+m_3^2+m_4^2}}{(q)_{2m_1}}
\qbins{m_1+m_2-m_3}{2m_2}\qbins{m_2+m_3-m_4}{2m_3}\qbins{m_3}{m_4} \\
=\frac{(q^{35},q^{36},q^{71};q^{71})_{\infty}}{(q;q)_{\infty}}
\end{multline*}
\and
\begin{equation*}
\sum_{n_1,n_2,m_3,m_4\geq 0}
\frac{q^{(n_1+n_2+m_3)^2+(n_2+m_3)^2+m_3^2+m_4^2}}
{(q)_{n_1}(q)_{n_2}(q)_{2m_3}}\qbins{m_3}{m_4}
=\frac{(q^{14},q^{15},q^{29};q^{29})_{\infty}}{(q;q)_{\infty}}
\end{equation*}
to be compared with
\begin{equation*}
\sum_{n_1,\dots,n_{34}\geq 0}
\frac{q^{N_1^2+\cdots+N_{34}^2}}
{(q)_{n_1}\cdots(q)_{n_{34}}}
=\frac{(q^{35},q^{36},q^{71};q^{71})_{\infty}}{(q;q)_{\infty}}
\end{equation*}
\and
\begin{equation*}
\sum_{n_1,\dots,n_{14}\geq 0}
\frac{q^{N_1^2+\cdots+N_{14}^2}}
{(q)_{n_1}\cdots(q)_{n_{14}}}
=\frac{(q^{14},q^{15},q^{29};q^{29})_{\infty}}{(q;q)_{\infty}}.
\end{equation*}
The most efficient sum sides of course occur in the identities
involving the Fibonacci numbers, with
$\log(2ab+1)/d(a,b)=\log(2F_kF_{k-1}+1)/(k-2)\to 2\log((1+\sqrt{5})/2)$
when $k$ tends to infinity.

We have not yet come to the end of our list of corollaries to
Theorem~\ref{thmmain} and next we are going to exploit the fact
that the polynomials in \eqref{eqmain} are not reciprocal.
If we define the polynomials $f_{a,b}(L,M)$ exactly as in \eqref{Fab}
and \eqref{Faab} but change the term $mC(a,b)m$ to $\bar{m}C(a,b)m$ with 
$\bar{m}=(m_1,\dots,m_{d(a,b)-1},0)$,
i.e., $\bar{m}C(a,b)m=mC(a,b)m+m_{d(a,b)}(m_{d(a,b)-1}-m_{d(a,b)})$
and define the special case $f_{2,1}(L,M)$ as
\begin{equation*}
f_{2,1}(L,M)=\sum_{m\geq 0}q^{L m}\qbin{2L+M-m}{2L}\qbin{L}{m},
\end{equation*}
then the $q\to 1/q$ version of Theorem~\ref{thmmain} can be states as follows.
\begin{corollary}\label{mainrecip}
For $L,M$ nonnegative integers and $(a,b)$ a pair of coprime integers
such that $1\leq b<a$ there holds
\begin{equation*}
\sum_{j\in\Z}(-1)^j q^{\frac{1}{2}j((2ab-1)j+1)}\B(L,M,aj,bj)=f_{a,b}(L,M).
\end{equation*}
\end{corollary}
The proof is obvious and merely requires \eqref{qbininv} and \eqref{Brecip}
and the observation that $\bar{m}C(a,b)m=\sum_{j,k=1}^{d(a,b)}(\tau_j-1)m_j
C(a,b)_{j,k}m_k$.

The case $(a,b)=(3,1)$ deserves special attention. Taking $\cf(3,1)=[1,1]$
and replacing $m_1\to L-i-n$ and $m_2 \to i$ gives
\begin{equation*}
\sum_{j\in\Z}(-1)^j q^{\frac{1}{2}j(5j+1)}\B(L,M,3j,j)=
\sum_{n,i\geq 0} q^{n^2+i(L+n)}\qbins{2L+M-n-i}{2L}
\qbins{2L-2i-n}{n}\qbins{L-i-n}{i}.
\end{equation*}
This is our second doubly-bounded analogue of the first Rogers--Ramanujan
identity, reducing to \eqref{RR1} in the large $M$ limit.
More generally, when $(a,b)=(k+1,1)$ we obtain new doubly-bounded analogues
of the Andrews--Gordon identity \eqref{AGid}
\begin{multline*}
\sum_{j\in\Z}(-1)^j q^{\frac{1}{2}j((2k+1)j+1)}\B(L,M,(k+1)j,j)=
\sum_{n_1,\dots,n_{k-1},i\geq 0} q^{N_1^2+\cdots+N_{k-1}^2+
i(L+\tilde{N}_k)} \\
\times
\qbins{2L+M-N_1-i}{2L}\qbins{L-i(k-1)-\tilde{N}_k}{i}
\prod_{j=1}^{k-1} \qbins{2L-2ij-N_j-N_{j+1}-2\tilde{N}_j}{n_j},
\end{multline*}
with $N_k=0$ and $\tilde{N}_j=N_1+\cdots+N_{j-1}=
\sum_{l=1}^{k-1}\min(j-1,l)n_l$.

We now proceed exactly as before.
First, defining $f_{a,b}(L)$ as in \eqref{FabL} and \eqref{FaabL} 
but with $mC(a,b)m$ replaced by $\bar{m}C(a,b)m$ (and 
$f_{2,1}(L)$ as the right side of \eqref{GB1}) 
we find the following large $M$ limit result.
\begin{corollary}
For $L$ a nonnegative integer and $(a,b)$ a pair of coprime integers such
that $1\leq b<a$ there holds
\begin{equation}\label{Gf}
G(L,L;b-1/a,b,a)=
\sum_{j\in\Z}(-1)^j q^{\frac{1}{2}j((2ab-1)j+1)}\qbin{2L}{L-aj}=f_{a,b}(L).
\end{equation}
\end{corollary}
A different route to this corollary is to take \eqref{GF}, replace $q$ by
$1/q$ using \eqref{Ginv} and \eqref{qbininv}, and to then replace $b$ by
$a-b$. From this remark it is clear that \eqref{Gf} does not yield new 
positivity results independent of \eqref{corGenBor}.

To take the large $L$ limit in Corollary~\ref{mainrecip}
we define $\tilde{f}_{a,b}(M)$ as in \eqref{FaabM} and \eqref{FabM}
but with the usual $m_{d(a,b)}(m_{d(a,b)-1}-m_{d(a,b)})$ added to
the exponent of $q$.
There is one exception, namely, $\tilde{f}_{a,1}(M)=\tilde{F}_{a-1,1}(M)$
(with $\tilde{f}_{2,1}(M)=(q)_{2M}/(q)_M$).
\begin{corollary}
For $M$ a nonnegative integer and $(a,b)$ a pair of coprime integers such 
that $1\leq b<a$ there holds
\begin{equation*}
G(M,M;a-1/b,a,b)=\sum_{j\in\Z}(-1)^j q^{\frac{1}{2}j((2ab-1)j+1)}
\qbin{2M}{M-bj}=\tilde{f}_{a,b}(M).
\end{equation*}
\end{corollary}
Letting $M$ tend to infinity we have reached our last theorem of this section.
\begin{theorem}\label{RRab2}
For $|q|<1$ and $(a,b)$ a pair of coprime integers such that $1\leq b<a$
there holds
\begin{equation*}
\sum_{m\in\Z^{d(a,b)}_+}\frac{q^{\bar{m}C(a,b)m}}{(q)_{2m_1}}
\prod_{j=2}^{d(a,b)}\qbins{\tau_j m_j+n_j}{\tau_j m_j}
=\frac{(q^{ab-1},q^{ab},q^{2ab-1};q^{2ab-1})_{\infty}}{(q;q)_{\infty}}
\end{equation*}
for $a<2b$, and
\begin{multline*}
\sum_{\substack{n_1,\dots,n_{a_0}\geq 0 \\
m_{a_0+1},\dots,m_{d(a,b)}\geq 0}}
\frac{q^{(N_1+m_{a_0+1})^2+\cdots+(N_{a_0}+m_{a_0+1})^2}}
{(q)_{n_1}\cdots(q)_{n_{a_0}}(q)_{2m_{a_0+1}}} \\
\times q^{\sum_{j,k=a_0+1}^{d(a,b)}
\bar{m}_j C(a,b)_{j,k}m_k}
\prod_{j=a_0+2}^{d(a,b)}\qbins{\tau_j m_j+n_j}{\tau_j m_j}
=\frac{(q^{ab-1},q^{ab},q^{2ab-1};q^{2ab-1})_{\infty}}{(q;q)_{\infty}}
\end{multline*}
for $a>2b\neq 2$ with $\bar{m}_j=m_j(1-\delta_{j,d(a,b)})$.
\end{theorem}
As examples one finds Theorem~\ref{Fib2} when
$(a,b)=(F_k,F_{k-1})$, Theorem~\ref{Fib4} when $(a,b)=(F_k,F_{k-2})$,
\begin{multline*}
\sum_{m_1,\dots,m_{k-1}\geq 0}
\frac{q^{M_1^2+\cdots+M_{k-2}^2-m_{k-2}M_{k-1}}}{(q)_{2m_1}}
\Bigl(\prod_{j=2}^{k-2}\qbins{m_{j-1}+m_{j+1}}{2m_j}\Bigr)
\qbins{m_{k-2}}{m_{k-1}} \\
=\frac{(q^{k(k-1)-1},q^{k(k-1)},q^{2k(k-1)-1};q^{2k(k-1)-1})_{\infty}}
{(q;q)_{\infty}}
\end{multline*}
with $M_j=m_j-m_{j-1}$ ($m_0=0$) when $(a,b)=(k,k-1)$ ($k\geq 3)$,
\begin{multline*}
\sum_{m_1,\dots,m_4\geq 0}
\frac{q^{m_1^2+(m_1-m_2)^2+m_3^2+m_3 m_4}}{(q)_{2m_1}}
\qbins{m_1+m_2-m_3}{2m_2}\qbins{m_2+m_3-m_4}{2m_3}\qbins{m_3}{m_4} \\
=\frac{(q^{34},q^{35},q^{69};q^{69})_{\infty}}{(q;q)_{\infty}}
\end{multline*}
when $(a,b)=(7,5)$, and
\begin{equation*}
\sum_{n_1,n_2,m_3,m_4\geq 0}
\frac{q^{(n_1+n_2+m_3)^2+(n_2+m_3)^2+m_3^2+m_3 m_4}}
{(q)_{n_1}(q)_{n_2}(q)_{2m_3}}\qbins{m_3}{m_4}
=\frac{(q^{13},q^{14},q^{27};q^{27})_{\infty}}{(q;q)_{\infty}}
\end{equation*}
when $(a,b)=(7,2)$.

\section{Proof of Theorem \ref{thmmain}}\label{secproof}
Throughout we assume that $(a,b)$ is a pair of positive, coprime integers
such that $a>b$. When $(a,b)\neq (2,1)$ we will, for definiteness, choose 
the representation of $\cf(a,b)$ with $a_n\geq 2$.
Now recall the definition of $d(a,b)$ as the sum of the partial
quotients of $\cf(a,b)$ and note that $d(2,1)=1$ and $d(a,b)>1$ 
if $(a,b)\neq (2,1)$.
Also note that for $(a,b)=(2,1)$ the theorem is nothing but
the identity \eqref{Bnewp}. (To see this recall that we
already calculated $F_{2,1}(L,M)$ and that it coincides with the
right side of \eqref{Bnewp}.)
We may use these facts to set up a proof by induction on $d(a,b)$.

Let us now assume that the theorem is valid for all $(a',b')$ such that
$d(a',b')=d$ and use this to prove its validity for all $(a,b)$ such that 
$d(a,b)=d+1$. There are four different cases to be considered depending
on the relative values of $a$ and $b$. Since we already dealt with 
$(a,b)=(2,1)$ we may assume that $a\neq 2b$. For brevity let us denote the
left-hand side of \eqref{eqmain} by $B_{a,b}(L,M)$.
\subsection{Proof for $a<2b$}
Let $(a,b)$ be a pair such that $a<2b$ and $d(a,b)=d+1$. Since $a<2b$ it
follows from \eqref{B2} that $B_{a,b}(L,M)$ satisfies the recurrence
\begin{equation}\label{Brec}
B_{a,b}(L,M)=\sum_{i\geq 0}q^{i^2}\qbin{2L+M-i}{2L}B_{b,a-b}(i,L-i).
\end{equation}
\subsubsection{Proof for $\frac{3}{2}b<a<2b$}\label{ss1}
If we write $a'=b$ and $b'=a-b$ then the condition $3b<2a$ translates
into $a'<2b'$ and hence $\cf(a',b')=b'/(a'-b')$ which we denote by
$[a'_0,\dots,a_{n'}']$.  Since $\cf(a,b)=b/(a-b)=1+1/\cf(a',b')$
we conclude that $\cf(a,b)=[1,a'_0,\dots,a_{n'}']$ and $d(a',b')=
d(a,b)-1=d$. (Note that by excluding $(a,b)=(3,2)$, which gives 
$(a',b')=(2,1)$, we avoid complications due to the fact that,
by our choice of representation, 
$\cf(3,2)=[2]=[1+a'_0]$ and not $\cf(3,2)=[1,1]=[1,a'_0]$.)
Consequently we may use our induction hypothesis to
replace $B_{a',b'}$ in \eqref{Brec} by $F_{a',b'}$
defined in \eqref{Fab}.
Abbreviating $\I(a',b')$ and $C(a',b')$ by $\I'$ and $C'$ this yields
\begin{equation*}
B_{a,b}(L,M)=\sum_{i\geq 0}q^{i^2}\qbins{2L+M-i}{2L}
\sum_{m'\in\Z_+^d}q^{m'C'm'}\qbins{L+i-m'_1}{2i}
\prod_{j=1}^d\qbins{\tau_j m'_j+n'_j}{\tau_j m'_j},
\end{equation*}
with
\begin{equation}\label{nprime}
n'_j=i\delta_{j,1}-\sum_{k=1}^d C'_{j,k} m'_k
\quad \text{for $j=1,\dots,d$.}
\end{equation}
Next rename $i$ as $m_1$ and $m'_j$ as $m_{j+1}$, and define
\begin{equation}\label{IIp}
\I=\left(\begin{array}{c|ccc}1&-1&& \\ \hline 1&&& \\ &&\I'& \\ &&& \\
\end{array}\right)
\end{equation}
and $C=2I-\I$ with $I$ the $(d+1)\times (d+1)$ identity matrix.
Observing that $\I=\I(a,b)$ and $C=C(a,b)$ thanks to 
$\cf(a,b)=1+1/\cf[a',b']=[1,a_0',\dots,a_{n'}']$, this gives
\begin{equation}\label{Bprime}
B_{a,b}(L,M)=\sum_{m\in\Z_+^{d+1}}q^{mC(a,b)m}\qbins{2L+M-m_1}{2L}
\prod_{j=1}^{d+1}\qbins{\tau_j m_j+n_j}{\tau_j m_j},
\end{equation}
with 
\begin{equation}\label{ndp1}
n_j=L\delta_{j,1}-\sum_{k=1}^{d+1} C_{j,k}(a,b) m_k
\quad \text{for $j=1,\dots,d+1$.}
\end{equation}
Since the right-hand side is exactly expression \eqref{Fab} for $F_{a,b}(L,M)$
with $d(a,b)=d+1$, this establishes \eqref{eqmain} for
$3b/2<a<2b$ and $d(a,b)=d+1$.

\subsubsection{Proof for $a\leq \frac{3}{2}b$}\label{ss2}
Again we write $a'=b$ and $b'=a-b$ but this time the condition $2a\leq 3b$
yields $a'\geq 2b'$. Hence $\cf(a',b')=(a'-b')/b'$ which is denoted by
$[a'_0,\dots,a_{n'}']$.  Since $\cf(a,b)=b/(a-b)=1+\cf(a',b')$
we conclude that $\cf(a,b)=[a'_0+1,\dots,a_{n'}']$ and $d(a',b')=
d(a,b)-1=d$. We may thus use the induction hypothesis to
replace $B_{a',b'}$ in \eqref{Brec} by $F_{a',b'}$ defined in \eqref{Faab}.
Abbreviating $\I(a',b')$ and $C(a',b')$ by $\I'$ and $C'$ this yields
\begin{equation*}
B_{a,b}(L,M)=\sum_{i\geq 0}q^{i^2}\qbins{2L+M-i}{2L}
\sum_{m'\in\Z_+^d}q^{i(i-2m_1)+m'C'm'}\qbins{L+m'_1}{2i}
\prod_{j=1}^d\qbins{\tau_j m'_j+n'_j}{\tau_j m'_j},
\end{equation*}
with $n'_j$ again given by \eqref{nprime}.
Once more we change $i\to m_1$ and $m'_j\to m_{j+1}$, and define
\begin{equation}\label{IIp2}
\I=\left(\begin{array}{c|ccc}0&1&& \\ \hline 1&&& \\ &&\I'& \\ &&& \\
\end{array}\right)
\end{equation}
and $C=2I-\I$. It follows from $\cf(a,b)=1+\cf[a',b']=[a_0'+1,\dots,a_{n'}']$
that $\I=\I(a,b)$ and $C=C(a,b)$, again leading to 
\eqref{Bprime} and \eqref{ndp1}, thus establishing \eqref{eqmain}
for $a\leq 3b/2$ and $d(a,b)=d+1$.

\subsection{Proof for $a>2b$}
Let $(a,b)$ be an admissible pair such that $a>2b$ and $d(a,b)=d+1$.
Since $a>2b$ it follows from \eqref{B1} that $B_{a,b}(L,M)$ satisfies
the recurrence
\begin{equation}\label{Brec2}
B_{a,b}(L,M)=\sum_{i\geq 0}q^{i^2}\qbin{2L+M-i}{2L}B_{a-b,b}(L-i,i).
\end{equation}
\subsubsection{Proof for $2b<a<3b$}
If we write $a'=a-b$ and $b'=b$ then the condition $a<3b$ implies that
$a'<2b'$, and we may copy the first part of the first paragraph
of section~\ref{ss1} and replace $B_{a',b'}$ in \eqref{Brec2} by 
$F_{a',b'}$ defined in \eqref{Fab}. With the same notation as before
this leads to
\begin{equation*}
B_{a,b}(L,M)=\sum_{i\geq 0}q^{i^2}\qbins{2L+M-i}{2L}
\sum_{m'\in\Z_+^d}q^{m'C'm'}\qbins{2L-i-m'_1}{2L-2i}
\prod_{j=1}^d\qbins{\tau_j m'_j+n'_j}{\tau_j m'_j},
\end{equation*}
with
\begin{equation}\label{nprime2}
n'_j=(L-i)\delta_{j,1}-\sum_{k=1}^d C'_{j,k} m'_k
\quad \text{for $j=1,\dots,d$.}
\end{equation}
This time we relabel $i$ as $L-m_1$ and $m'_j$ as $m_{j+1}$, and define $\I$
as in \eqref{IIp} and $C=2I-\I$. Since $\I=\I(a,b)$ and $C=C(a,b)$ thanks to 
$\cf(a,b)=1+1/\cf[a',b']$ this gives
\begin{equation}\label{Bprime2}
B_{a,b}(L,M)=\sum_{m\in\Z_+^{d+1}}q^{L(L-2m_1)+mC(a,b)m}\qbins{L+M+m_1}{2L}
\prod_{j=1}^{d+1}\qbins{\tau_j m_j+n_j}{\tau_j m_j},
\end{equation}
with $n_j$ given by \eqref{ndp1}.
Since the right-hand side is exactly expression \eqref{Faab} for $d(a,b)=d+1$,
this results in \eqref{eqmain} for $2b<a<3b$ and $d(a,b)=d+1$.
\subsubsection{Proof for $a\geq 3b$}
Writing $a'=a-b$ and $b'=b$ the condition $a\geq 3b$ becomes $a'\geq 2b'$,
and we may copy the first part of the first paragraph of section~\ref{ss2}
and replace $B_{a',b'}$ in \eqref{Brec2} by $F_{a',b'}$ defined in
\eqref{Faab}. With the same notation as before this leads to
\begin{equation*}
B_{a,b}(L,M)=\sum_{i\geq 0}q^{i^2}\qbins{2L+M-i}{2L}
\sum_{m'\in\Z_+^d}q^{m'C'm'}\qbins{L+m'_1}{2L-2i}
\prod_{j=1}^d\qbins{\tau_j m'_j+n'_j}{\tau_j m'_j},
\end{equation*}
with $n_j'$ given by \eqref{nprime2}.
Renaming $i\to L-m_1$ and $m'_j \to m_{j+1}$, defining $\I$
as in \eqref{IIp} and $C=2I-\I$, and observing that $\I=\I(a,b)$ and 
$C=C(a,b)$ because $\cf(a,b)=1+\cf[a',b']$, this again gives
\eqref{Bprime2} and \eqref{ndp1} establishing \eqref{eqmain}
for $a\geq 3b$ and $d(a,b)=d+1$ and completing the proof.

\section{Further positivity results}
So far we have no identities for $G(N,M;\alpha,\beta,K)$ 
with both $\alpha$ and $\beta$ noninteger. To obtain such
results our starting point will be yet another doubly-bounded analogue
of the first Rogers--Ramanujan identity.
\begin{lemma}
For $L,M\geq 0$ there holds
\begin{equation}\label{Bnewp2}
\sum_{j\in\Z}(-1)^j q^{\frac{1}{2}j(5j+1)}\B(L,M,2j+1,j)
=\sum_{n\geq 0}q^{n^2}\qbins{2L+M-n-1}{2L-1}\qbins{L-1}{n}.
\end{equation}
\end{lemma}

\begin{proof}
As a first step we add zero in the form
\begin{equation*}
q^{L+1}\sum_{j\in\Z}(-1)^j q^{\frac{5}{2}j(j+1)}\qbins{L+M+j}{M-j-1}
\qbins{L+M-j-1}{M+j}
\end{equation*}
to the left side of \eqref{Bnewp2}.
By the recurrence \eqref{qbinrec} we are then to prove
\begin{equation}\label{Brep}
\sum_{j\in\Z}(-1)^j q^{\frac{1}{2}j(5j+1)}\qbins{L+M+j}{M-j}
\qbins{L+M-j-1}{M+j}
=\sum_{n\geq 0}q^{n^2}\qbins{2L+M-n-1}{2L-1}\qbins{L-1}{n}.
\end{equation}
Suppressing the $L$-dependence in the first, we denote the left sides of
\eqref{Brep} and \eqref{Bnewp} by $f_M$ and $g_{L,M}$, respectively.
Using \eqref{qbinrec} we then get
\begin{align*}
f_M&=g_{L-1,M}+
q^L\sum_{j\in\Z}(-1)^j q^{\frac{5}{2}j(j+1)}\qbins{L+M+j-1}{M-j-1}
\qbins{L+M-j-1}{M+j} \\
&=g_{L-1,M}+q^{2L-1}f_{M-1}
+q^L\sum_{j\in\Z}(-1)^j q^{\frac{5}{2}j(j+1)}\qbins{L+M+j-1}{M-j-1}
\qbins{L+M-j-2}{M+j} \\
&=g_{L-1,M}+q^{2L-1}f_{M-1}.
\end{align*}
Next we let $f_M$ and $g_{L,M}$ denote the right sides of \eqref{Brep}
and \eqref{Bnewp}. One application of \eqref{qbinrec} shows that
the same recurrence again holds.
Since \eqref{Brep} is true for $M=0$ we are done.
\end{proof}

Now that \eqref{Bnewp2} has been proven we closely follow the
work of the previous two sections.
The present situation is, however, notationally more involved
and further definitions related to continued fractions are needed.
Let $(a,b)$ be the usual ordered pair of coprime integers
with associated continued fraction $\cf(a,b)=(a/b-1)^{\sign(a-2b)}
=[a_0,\dots,a_n]$. In principle we could still allow for both representations
of $\cf(a,b)$ but many of the equations below are sensitive to the chosen
representation and to avoid unnecessary complications we demand that 
$a_n\geq 2$ for $(a,b)\neq (2,1)$. Given $(a,b)$ we define a second pair 
$(\bar{a},\bar{b})$ of positive, coprime integers as follows
\begin{equation}\label{bar}
\frac{\bar{a}}{\bar{b}}=\begin{cases}
[1,a_0,\dots,a_{n-1}]=1+1/[a_0,\dots,a_{n-1}] & \text{for $a<2b$} \\
 {} [a_0+1,a_1,\dots,a_{n-1}]=1+[a_0,\dots,a_{n-1}] & \text{for $a>2b$,}
\end{cases}
\end{equation}
with special cases $(\bar{a},\bar{b})=(1,0)$ for $(a,b)=(a,1)$ ($a\geq 2$)
and $(\bar{a},\bar{b})=(1,1)$ for $(a,b)=(a,a-1)$ ($a>2$).
Since $1/[1,c_0,\dots,c_n]+1/[c_0+1,\dots,c_n]=1$ it readily follows
that for $(a_1,b_1)$ and $(a_2,b_2)$ such that $a_1=a_2$ and $b_2=a_1-b_1$
there holds $\bar{a}_1=\bar{a}_2=\bar{b}_1+\bar{b}_2$.
It is also easy to see that if $a<2b$ then $\bar{a}\leq 2\bar{b}$
and if $a\geq 2b$ then $\bar{a}\geq 2\bar{b}$.
We also observe that $[1,a_0,\dots,a_n]=a/b$ ($a\geq 2b$) and
$[a_0+1,a_1,\dots,a_n]=a/b$ ($a\leq 2b$). 
In the language of continued fractions this means that $\bar{a}/\bar{b}$
is the $n$th convergent of the continued fraction of $a/b$ (which itself
is of order $n+1$). Care should however be taken with the anomalous case 
$(a,b)=(a,1)$ for which $a/b=a=[a]$ and $(\bar{a},\bar{b})=(1,0)$.
As an example of the above definitions let $(a,b)=(19,12)$.
Then $\cf(19,12)=[1,1,2,2]$, $\bar{a}/\bar{b}=[1,1,1,2]=8/5$ and 
$a/b=[1,1,1,2,2]$. Similarly, if $(a,b)=(19,7)$, then $\cf(19,7)=[1,1,2,2]$,
$\bar{a}/\bar{b}=[2,1,2]=8/3$ and $a/b=[2,1,2,2]$.

We further need to define the analogues of the polynomials $F_{a,b}(L,M)$,
which will be denoted by $H_{a,b}(L,M)$.
For $(a,b)=(2,1)$ we take it to be the right-hand side of \eqref{Bnewp2}
and for all other $(a,b)$
\begin{equation}\label{Hab}
H_{a,b}(L,M)=\sum_{m\in\Z_+^{d(a,b)}}q^{mC(a,b)m+A_m}\qbins{2L+M-m_1}{2L}
\prod_{j=1}^{d(a,b)}\qbins{\tau_j m_j+n_j-\delta_{j,d(a,b)}}
{\tau_j m_j-\delta_{j,d(a,b)-1}}
\end{equation}
for $a<2b$, and
\begin{equation}\label{Haab}
H_{a,b}(L,M)=\sum_{m\in \Z_+^{d(a,b)}}q^{L(L-2m_1)+mC(a,b)m+A_m}
\qbins{L+M+m_1}{2L}
\prod_{j=1}^{d(a,b)}\qbins{\tau_j m_j+n_j-\delta_{j,d(a,b)}}
{\tau_j m_j-\delta_{j,d(a,b)-1}}
\end{equation}
for $a>2b$. Here $mC(a,b)m$ and $n_j$ are still given by \eqref{mCm}
and \eqref{mn}. The term $A_m$ denotes the the linear term 
$2m_{d(a,b)}-2m_{d(a,b)-1}+1$.
Since $a_n\geq 2$ we thus have
\begin{equation*}
mC(a,b)m+A_m=\sum_{j=0}^n
\Bigl(m_{t_j+1}^2+\sum_{k=t_j+1}^{t_{j+1}-1}
(m_{k}-m_{k+1}-\delta_{k,d(a,b)-1})^2\Bigr).
\end{equation*}

The analogue of theorem~\ref{thmmain} now breaks up into two separate
statements.
\begin{theorem}\label{thmmain2}
For $L,M$ nonnegative integers and $(a,b)$ a pair of coprime integers
such that $1\leq b<a$ there holds (i)
\begin{equation}\label{Ex1}
\sum_{j\in\Z}(-1)^j q^{\frac{1}{2}j((2ab+1)j+4\bar{a}b+1)+\bar{a}\bar{b}}
\B(L,M,aj+\bar{a},bj+\bar{b})=H_{a,b}(L,M)
\end{equation}
for $a<2b$ and $\cf(a,b)$ a continued fraction of even order (i.e., $n$ even),
or $a>2b$ and $\cf(a,b)$ a continued fraction of odd order, (ii)
\begin{equation}\label{Ex2}
\sum_{j\in\Z}(-1)^j q^{\frac{1}{2}j((2ab+1)j+4a\bar{b}+1)+\bar{a}\bar{b}}
\B(L,M,aj+\bar{a},bj+\bar{b})=H_{a,b}(L,M)
\end{equation}
for $a<2b$ and $\cf(a,b)$ of odd order, or $a\geq 2b$ and $\cf(a,b)$ of
even order.
\end{theorem}
Instead of boring the reader with a complete list of
analogues of the results of section~\ref{secGB},
we restrict ourselves to the analogue of Corollary~\ref{corGenBor}.
This follows after letting $M$ tend to infinity in Theorem~\ref{thmmain2} 
and using \eqref{Blim}.
\begin{corollary}\label{cormain2}
Let $(a,b)$ be a pair of coprime integers such that $1\leq b<a$,
with $(a,b)=(2,1)$ excluded.
Then (i) $G(L+\bar{a},L-\bar{a};b-2\bar{a}b/a,b+1/a+2\bar{a}b/a,a)$ is a 
polynomial with nonnegative coefficients if $a<2b$ and $\cf(a,b)$ has 
even order, or $a>2b$ and $\cf(a,b)$ has odd order, (ii)
$G(L+\bar{a},L-\bar{a};b-2\bar{b},b+1/a+2\bar{b},a)$ is a 
polynomial with nonnegative coefficients if $a<2b$ and $\cf(a,b)$ has 
odd order, or $a>2b$ and $\cf(a,b)$ has even order.
\end{corollary}
No further positivity results arise from Theorem~\ref{thmmain2}.
If we replace $q$ by $1/q$ and then let $M$ tend to infinity we
get results which are equivalent to those obtained by exploiting the
duality \eqref{Ginv} in Corollary~\ref{cormain2}. Although this is
similar to the situation encountered in section~\ref{secGB} it is
more cumbersome to prove this. A brief derivation proceeds as follows.
Start with \eqref{Ex2} and replace $q\to 1/q$ using \eqref{Brecip}.
This gives the polynomial identity
\begin{equation*}
\sum_{j\in\Z}(-1)^j q^{\frac{1}{2}j((2ab-1)j+4a\bar{b}-1)+\bar{a}\bar{b}}
\B(L,M,aj+\bar{a},bj+\bar{b})=q^{2LM}H_{a,b}(L,M;1/q).
\end{equation*}
Letting $M$ tend to infinity establishes the positivity of
$G(L-\bar{a},L+\bar{a};b-1/a+2\bar{b},b-2\bar{b},a)$.
All of this is of course for $a<2b$ and $\cf(a,b)$ of even order,
or $a>2b$ and $\cf(a,b)$ of odd order.
Next replace $b\to a-b$ which has the effect of changing
$\bar{b}\to \bar{a}-\bar{b}$. Hence
$G(L-\bar{a},L+\bar{a};a-b-1/a+2(\bar{a}-\bar{b}),a-b-2(\bar{a}-\bar{b}),a)$
is positive for $a>2b$ and $\cf(a,b)$ of even order,
or $a<2b$ and $\cf(a,b)$ of odd order.
Letting $q\to 1/q$ in item (ii) of the Corollary~\ref{cormain2}
using \eqref{Ginv} yields the same result.
Starting with \eqref{Ex2} instead of \eqref{Ex1} reproduces
the $q\to 1/q$ analogue of item (i) of the corollary in much the same way.

Before we come to the proof of Theorem~\ref{thmmain2} let us consider
the example $(a,b)=(3,1)$. Then $(\bar{a},\bar{b})=(1,0)$ and $\cf(3,1)=[2]$
which is of order $0$. Since also $a>2b$ we are to use \eqref{Ex2}
leading to
\begin{multline}\label{RR2inv}
\sum_{j\in\Z}(-1)^j q^{\frac{1}{2}j(7j+1)}\B(L,M,3j+1,j) \\
=\sum_{m_1,m_2\geq 0}q^{(L-m_1)^2+(m_1-m_2-1)^2}
\qbins{L+M+m_1}{2L}\qbins{L+m_2}{2m_1-1}\qbins{m_1-1}{m_2}.
\end{multline}
The reason for giving this rather atypical example with trivial
$(\bar{a},\bar{b})$ is its $q\to 1/q$ counterpart, which after the 
replacements $m_1\to L-i-n$ and $m_2\to L-2i-n-1$ becomes
\begin{multline*}
\sum_{j\in\Z}(-1)^j q^{\frac{1}{2}j(5j+3)}\B(L,M,3j+1,j) \\
=\sum_{n,i\geq 0}q^{n(n+1)+i(L+n+1)}
\qbins{2L+M-i-n}{2L}\qbins{2L-2i-n-1}{n}\qbins{L-i-n-1}{i}.
\end{multline*}
This doubly-bounded analogue of the second Rogers--Ramanujan
yields \eqref{RR2} of the introduction in the large $M$ limit.

\section{Proof of Theorem \ref{thmmain2}}
Though considerably more involved, the proof proceeds along the
same lines as the proof of Theorem~\ref{thmmain} given
in section~\ref{secproof}. Again we carry out induction on $d(a,b)$,
but the first difference is that the case $(a,b)=(2,1)$ which has $d(2,1)=1$,
is special and is not included in either \eqref{Hab} or \eqref{Haab}.
We therefore have to first prove the cases $(a,b)=(3,1)$ and $(3,2)$
which are the only solutions to $d(a,b)=2$. These can then serve as
starting point for our induction. 
Applying \eqref{B2} to \eqref{Bnewp2} gives
\begin{multline}\label{ab32}
\sum_{j\in\Z}(-1)^j q^{\frac{1}{2}j(13j+9)+1}\B(L,M,3j+1,2j+1) \\
=\sum_{i,n\geq 0}q^{i^2+n^2}\qbins{2L+M-i}{2L}
\qbins{L+i-n-1}{2i-1}\qbins{i-1}{n}.
\end{multline}
For $(a,b)=(3,2)$ one finds $\cf(a,b)=[2]$ of even order, $a<2b$
and $(\bar{a},\bar{b})=(1,1)$. 
Hence the left side of \eqref{Ex1} for $(a,b)=(3,2)$ agrees with
the left side of \eqref{ab32}.
To see that also the right sides agree we rewrite the right side
of \eqref{ab32} in terms of the summation variables $m_1=i$ and 
$m_2=i-n-1$. This gives
\begin{equation*}
\sum_{m_1,m_2\geq 0}q^{m_1^2+(m_1-m_2-1)^2}\qbins{2L+M-m_1}{2L}
\qbins{L+m_2}{2m_1-1}\qbins{m_1-1}{m_2}
\end{equation*}
in accordance with the $(a,b)=(3,2)$ case of \eqref{Hab}.
Similarly, applying \eqref{B1} to \eqref{Bnewp2} gives
\begin{equation*}
\sum_{j\in\Z}(-1)^j q^{\frac{1}{2}j(7j+1)}\B(L,M,3j+1,j)
=\sum_{i,n\geq 0}q^{i^2+n^2}\qbins{2L+M-i}{2L}
\qbins{2L-i-n-1}{2L-2i-1}\qbins{L-i-1}{n}.
\end{equation*}
Introducing 
$m_1=L-i$ and $m_2=L-i-n-1$ on the right transforms this into 
\eqref{RR2inv} which is the $(a,b)=(3,1)$ case of Theorem~\ref{thmmain2}.

Now we are prepared for the induction step, and we assume the theorem 
to be correct for all $(a',b')$ such that $d(a',b')=d$ in order
to prove its validity for all $(a,b)$ such that $d(a,b)=d+1$. 
There are eight cases to be considered depending
on the relative values of $a$ and $b$ and on the order of $\cf(a,b)$.
For convenience let us introduce the notation
$B_{a,b}^{\text{e},<}(L,M)$ for the left side of \eqref{Ex1} when $a<2b$ and
$\cf(a,b)$ has even order. In the same way we define
$B_{a,b}^{\text{o},>}(L,M)$, $B_{a,b}^{\text{e},>}(L,M)$
and $B_{a,b}^{\text{o},<}(L,M)$.
\subsection{Proof for $\frac{3}{2}b<a<2b$ with $\cf(a,b)$ of even order}
Let $(a,b)$ be a pair such that $3b/2<a<2b$, $\cf(a,b)$ of even order $n$
and $d(a,b)=d+1$. For such $(a,b)$ we will show that
\begin{equation}\label{Brec1}
B_{a,b}^{\text{e},<}(L,M)=\sum_{i\geq 0}q^{i^2}\qbin{2L+M-i}{2L}
B_{b,a-b}^{\text{o},<}(i,L-i),
\end{equation}
with $d(b,a-b)=d$. First write $a'=b$ and $b'=a-b$. Then $3b<2a$ gives
$a'<2b'$ (in accordance with \eqref{Brec1}) and hence 
$\cf(a',b')=[a'_0,\dots,a_{n'}']=b'/(a'-b')$.
Since $\cf(a,b)=b/(a-b)=1+1/\cf(a',b')$ one finds 
$\cf(a,b)=[1,a'_0,\dots,a_{n'}']$.
This yields $n'=n-1$ is odd and $d(a',b')=d(a,b)-1=d$.
We may now conclude that at least the labels in \eqref{Brec1}
are correct. 
Next we use the definition of $B_{a',b'}^{\text{o},<}$ and the Burge
transform \eqref{B2} to compute the right-hand side of \eqref{Brec1} as
\begin{multline*}
\sum_{j\in\Z}(-1)^j 
q^{\frac{1}{2}j((2a'(a'+b')+1)j+4a'(\bar{a}'+\bar{b}')+1)
+\bar{a}'(\bar{a}'+\bar{b}')} \\
\times \B(L,M,(a'+b')j+\bar{a}'+\bar{b}',a'j+\bar{a}')
\end{multline*}
Since $\bar{a}/\bar{b}=[1,a_0,\dots,a_{n-1}]=[1,1,a_0',\dots,a_{n'-1}]=
1+1/[1,a_0,\dots,a_{n'-1}]=(\bar{a}'+\bar{b}')/\bar{a}'$ we find that
$\bar{a}'=\bar{b}$ and $\bar{b}'=\bar{a}-\bar{b}$.
Also using $a'=b$ and $b'=a-b$ we can simplify the above
expression to
\begin{equation*}
\sum_{j\in\Z}(-1)^j 
q^{\frac{1}{2}j((2ab+1)j+4\bar{a}b+1)+\bar{a}\bar{b}} 
\B(L,M,aj+\bar{a},bj+\bar{b}).
\end{equation*}
This is precisely the left side of \eqref{Brec1} as we set out to prove.

The remaining part of the proof proceeds exactly as
the proof given in section~\ref{ss1}. In a few words,
we can use the induction hypothesis to replace the right side of
\eqref{Brec1} by $H_{a',b'}(i,L-i)$ given in \eqref{Hab}.
Making the necessary variable changes gives the expression
claimed by the theorem.

\subsection{Proof for $\frac{3}{2}b<a<2b$ with $\cf(a,b)$ of odd order}
Let $(a,b)$ be a pair such that $3b/2<a<2b$, $\cf(a,b)$ of odd order $n$
and $d(a,b)=d+1$. Then \eqref{Brec1} is replaced by 
\begin{equation}\label{Brec3}
B_{a,b}^{\text{o},<}(L,M)=\sum_{i\geq 0}q^{i^2}\qbin{2L+M-i}{2L}
B_{b,a-b}^{\text{e},<}(i,L-i),
\end{equation}
with $d(b,a-b)=d$.
Again we define $a'=b$ and $b'=a-b$. 
Copying the paragraph below \eqref{Brec1} replacing the one occurrence of
`odd' by `even', shows that the labels are again correct.
Computing the right-hand side of \eqref{Brec3} using
$\bar{a}'=\bar{b}$ and $\bar{b}'=\bar{a}-\bar{b}$ results in
\begin{equation*}
\sum_{j\in\Z}(-1)^j 
q^{\frac{1}{2}j((2ab+1)j+4a\bar{b}+1)+\bar{a}\bar{b}} 
\B(L,M,aj+\bar{a},bj+\bar{b})
\end{equation*}
in agreement with the left-hand side of \eqref{Brec3}.
The rest of the proof follows that of section~\ref{ss1}.

\subsection{Proof of the remaining six cases}
The proofs of the remaining cases are simple modifications of 
the previous two and are therefore omitted. 
For completeness we just state the six key-identities
\begin{align*}
B_{a,b}^{\text{p},<}(L,M)&=\sum_{i\geq 0}q^{i^2}\qbin{2L+M-i}{2L}
B_{b,a-b}^{\text{p},>}(i,L-i) \quad \text{for $2a<3b$} \\
\end{align*}
and
\begin{align*}
B_{a,b}^{\text{p},>}(L,M)&=\sum_{i\geq 0}q^{i^2}\qbin{2L+M-i}{2L}
B_{a-b,b}^{\bar{\text{p}},<}(L-i,i) \quad \text{for $2b<a<3b$} \\
B_{a,b}^{\text{p},>}(L,M)&=\sum_{i\geq 0}q^{i^2}\qbin{2L+M-i}{2L}
B_{a-b,b}^{\text{p},>}(L-i,i) \quad \text{for $a>3b$,}
\end{align*}
where $\text{p}=\text{e,o}$, $\bar{\text{e}}=\text{o}$ and
$\bar{\text{o}}=\text{e}$.

\section{Rogers--Ramanujan-type identities for even moduli}
Closely related to the Andrews--Gordon identities \eqref{AGid}
are Bressoud's identities for even moduli \cite{Bressoud80}
\begin{equation}\label{Bressoudid}
\sum_{n_1,\dots,n_{k-1}\geq 0}
\frac{q^{N_1^2+\cdots+N_{k-1}^2}}{(q)_{n_1}\cdots(q)_{n_{k-2}}
(q^2;q^2)_{n_{k-1}}}
=\frac{(q^k,q^k,q^{2k};q^{2k})_{\infty}}{(q;q)_{\infty}}
\end{equation}
for $k\geq 2$, $N_j=n_j+\cdots+n_{k-1}$ and $|q|<1$.
An obvious question is whether also these identities can be
embedded in an infinite tree of Rogers--Ramanujan-type identities.
The answer to this is `yes' and is already implicit in Burge's original
paper on his transform.
In \cite[page 217]{Burge93} Burge states the following identity
for $L,M\geq 0$:
\begin{equation*}
\sum_{j\in\Z}(-1)^j q^{j^2}\B(L,M,j,j)=\qbin{L+M}{M}_{q^2}.
\end{equation*}
This identity is very similar to \eqref{Bnew} and if we define $I_{a,b}(L,M)$
by \eqref{Fab} and \eqref{Faab} but with each $q$-binomial
coefficient dressed with a subscript $q^{\bar{\tau}_j}$ where 
$\bar{\tau}_j=3-\tau_j$, then the following theorem is immediate.
\begin{theorem}
For $L,M$ nonnegative integers and $(a,b)$ a pair of coprime integers
such that $1\leq b<a$ there holds
\begin{equation*}
\sum_{j\in\Z}(-1)^j q^{abj^2}\B(L,M,aj,bj)=I_{a,b}(L,M).
\end{equation*}
\end{theorem}
As far as the generalized Borwein conjecture goes this is not
interesting, but letting both $L$ and $M$ tend to infinity 
yields an even modulus version of Theorem~\ref{RRab}.
\begin{theorem}\label{thmeven}
For $|q|<1$ and $(a,b)$ a pair of coprime integers such that $1\leq b<a$
there holds
\begin{equation*}
\sum_{m\in\Z^{d(a,b)}_+}\frac{q^{mC(a,b)m}}
{(q)_{2m_1}}
\prod_{j=2}^{d(a,b)}\qbins{\tau_j m_j+n_j}{\tau_j m_j}_{q^{\bar{\tau}_j}}
=\frac{(q^{ab},q^{ab},q^{2ab};q^{2ab})_{\infty}}{(q;q)_{\infty}}
\end{equation*}
for $a<2b$, and
\begin{multline*}
\sum_{\substack{n_1,\dots,n_{a_0}\geq 0 \\
m_{a_0+1},\dots,m_{d(a,b)}\geq 0}}
\frac{q^{(N_1+m_{a_0+1})^2+\cdots+(N_{a_0}+m_{a_0+1})^2}}
{(q)_{n_1}\cdots(q)_{n_{a_0-1}}
(q^{\bar{\tau}_{a_0}};q^{\bar{\tau}_{a_0}})_{n_{a_0}}
(q^{\bar{\tau}_{a_0+1}};q^{\bar{\tau}_{a_0+1}})_{\tau_{a_0+1}m_{a_0+1}}} \\
\times q^{\sum_{j,k=a_0+1}^{d(a,b)} m_j C(a,b)_{j,k} m_k}
\prod_{j=a_0+2}^{d(a,b)}\qbins{\tau_j m_j+n_j}{\tau_j m_j}_{q^{\bar{\tau}_j}}
=\frac{(q^{ab},q^{ab},q^{2ab};q^{2ab})_{\infty}}{(q;q)_{\infty}}
\end{multline*}
for $a\geq 2b$.
\end{theorem}
For $(a,b)=(k,1)$ this is Bressoud's \eqref{Bressoudid}.
The most interesting other examples again feature the Fibonacci numbers,
and for $(a,b)=(F_k,F_{k-1})$ and $(a,b)=(F_k,F_{k-2})$ one finds
\begin{multline*}
\sum_{m_1,\dots,m_{k-2}\geq 0}
\frac{q^{m_1^2+\cdots+m_{k-2}^2}}{(q)_{2m_1}}
\Bigl(\prod_{j=2}^{k-3}\qbins{m_{j-1}+m_j-m_{j+1}}{2m_j}\Bigr)
\qbins{m_{k-3}}{m_{k-2}}_{q^2} \\
=\frac{(q^{F_kF_{k-1}},q^{F_kF_{k-1}},q^{2F_kF_{k-1}};
q^{2F_kF_{k-1}})_{\infty}}{(q;q)_{\infty}}
\end{multline*}
for $k\geq 4$, and
\begin{multline*}
\sum_{m_1,\dots,m_{k-2}\geq 0}
\frac{q^{(m_1+m_2)^2+m_2^2+\cdots+m_{k-2}^2}}
{(q)_{m_1}(q)_{2m_2}}
\Bigl(\prod_{j=3}^{k-3}\qbins{m_{j-1}+m_j-m_{j+1}}{2m_j}\Bigr)
\qbins{m_{k-3}}{m_{k-2}}_{q^2} \\
=\frac{(q^{F_kF_{k-2}},q^{F_kF_{k-2}},q^{2F_kF_{k-2}};
q^{2F_kF_{k-2}})_{\infty}}{(q;q)_{\infty}}
\end{multline*}
for $k\geq 5$, where we have replaced $n_1\to m_1$ in comparison 
with Theorem~\ref{thmeven}.
These two series can again be extended to all $k\geq 3$ by taking
$\sum_{m\geq 0} q^{m^2}/(q^2;q^2)_m$ as left-hand sides when $k=3$
and $\sum_{m_1,m_2\geq 0} q^{(m_1+m_2)^2+m_2^2}/(q)_{m_1}(q^2;q^2)_{m_2}$ 
as left-hand side of the second series when $k=4$.

\section{Outlook: how tractable is the Borwein conjecture?}\label{secdis}
Using the Burge transform we have established two types of positivity
results for $G(N,M;\alpha,\beta,K)$.
The first type has $N=M$ with either $\alpha$ or $\beta$ being an integer,
and the second type has $N\neq M$ where both $\alpha$ and $\beta$ can be
noninteger.
This might lead one to suspect that proving the positivity of coefficients
of $G(n,n;\alpha,\beta,K)$ for noninteger $\alpha$ and $\beta$
is perhaps a much more difficult problem.
The sceptical reader might even doubt that nice identities for such
$G$ exist in the first place, casting doubt on the
claim made in the introduction that our failure to prove the positivity
of $A_n(q)=G(n,n;4/3,5/3,3)$ is possibly just a practical and not a
fundamental problem.

Indeed we believe that the original Borwein conjecture is quite a bit
deeper than the positivity results proven in this paper. However, in
our subsequent papers on this topic we will introduce new types of
transformations that settle more complicated cases of Bressoud's
conjecture. Some appealing examples are the final entries of the following
two sequences of identities:
\begin{align*}
G(n,n;1/2,2/2,2)&=\sum_{m_1\geq 0} q^{m_1 n}\qbin{n}{m_1} \\
G(n,n;3/3,4/3,3)&=\sum_{m_1,m_2\geq 0} q^{(n-m_1)^2+(m_1-m_2)^2}
\qbin{n+m_2}{2m_1}\qbin{m_1}{m_2} \\
G(n,n;5/4,6/4,4)&=\sum_{m_1,m_2,m_3\geq 0} 
q^{n(n-m_1)+m_2(m_2+m_3)}\qbin{n}{m_1}\qbin{m_1}{2m_2}\qbin{m_2}{m_3}
\end{align*}
and
\begin{align*}
G(n,n;2/2,3/2,2)&=\sum_{m_1\geq 0} q^{m_1^2}\qbin{n}{m_1} \\[1mm]
G(n,n;4/3,5/3,3)&=A_n(q)=\text{???} \\
G(n,n;6/4,7/4,4)&=\sum_{m_1,m_2,m_3\geq 0} 
q^{n(n-m_1)+m_2^2+m_3^2}\qbin{n}{m_1}\qbin{m_1}{2m_2}\qbin{m_2}{m_3}.
\end{align*}
The pattern is of course clear and we leave it to the reader to
fill in the missing item.

\subsection*{Note}
Alexander Berkovich kindly informed me that he has independently 
obtained several of the results established in this paper. 
In particular he has obtained the corollaries~\ref{Corab} and \ref{corGenBor}
and all results implied by these.

\bibliographystyle{amsplain}

\begin{thebibliography}{99}

\bibitem{Andrews74}
G. E. Andrews,
\textit{An analytic generalization of the Rogers--Ramanujan identities
for odd moduli},
Prod. Nat. Acad. Sci. USA \textbf{71} (1974), 4082--4085.

\bibitem{Andrews76}
G. E. Andrews,
\textit{The Theory of Partitions},
Encyclopedia of Mathematics and its Applications, Vol.~2,
(Addison-Wesley, Reading, Massachusetts, 1976).

\bibitem{Andrews95}
G. E. Andrews,
\textit{On a conjecture of Peter Borwein},
J. Symbolic Comput. \textbf{20} (1995), 487--501.

\bibitem{ABBBFV87}
G. E. Andrews, R. J. Baxter, D. M. Bressoud, W. H. Burge, P. J. Forrester
and G. Viennot,
\textit{Partitions with prescribed hook differences},
Europ. J. Combinatorics \textbf{8} (1987), 341--350.

\bibitem{BM96}
A. Berkovich and B. M. McCoy,
\textit{Continued fractions and fermionic representations
for characters of $M(p,p')$ minimal models},
Lett. Math. Phys. \textbf{37} (1996), 49--66.

\bibitem{BMS98}
A. Berkovich, B. M. McCoy and A. Schilling,
\textit{Rogers-Schur-Ramanujan type identities for the $M(p,p')$ minimal
models of conformal field theory},
Commun. Math. Phys. \textbf{\bf 191} (1998), 325--395.

\bibitem{Bressoud80}
D. M. Bressoud,
\textit{An analytic generalization of the Rogers--Ramanujan identities
with interpretation},
Quart. J. Maths. Oxford (2) \textbf{31} (1980), 385--399.

\bibitem{Bressoud81}
D. M. Bressoud,
\textit{Some identities for terminating $q$-series},
Math. Proc. Camb. Phil. Soc. \textbf{89} (1981), 211--223.

\bibitem{Bressoud96}
D. M. Bressoud,
\textit{The Borwein conjecture and partitions with prescribed hook
differences},
Electron. J. Combin. {\bf 3} (1996), \#4.

\bibitem{Burge93}
W. H. Burge,
\textit{Restricted partition pairs},
J. Combin. Theory Ser. A \textbf{63} (1993), 210--222.

\bibitem{FLW97}
O. Foda, K. S. M. Lee and T. A. Welsh,
\textit{A Burge tree of Virasoro-type polynomial identities},
Int. J. Mod. Phys. A \textbf{13} (1998), 4967--5012.

\bibitem{GR90}
G. Gasper and M. Rahman,
\textit{Basic Hypergeometric Series},
Encyclopedia of Mathematics and its Applications, Vol.~35,
(Cambridge University Press, Cambridge, 1990).

\bibitem{IKS99}
M. E. H. Ismail, D. Kim and D. Stanton,
\textit{Lattice paths and positive trigonometric sums},
Constr. Approx. \textbf{15} (1999), 69--81.

\bibitem{Paule85}
P. Paule,
\textit{On identities of the Rogers--Ramanujan type},
J. Math. Anal. Appl. \textbf{107} (1985), 255--284.

\bibitem{Rogers94}
L. J. Rogers,
\textit{Second memoir on the expansion of certain infinite products},
Proc. London Math. Soc. \textbf{25} (1894), 318--343.

\bibitem{Rogers17}
L. J. Rogers,
\textit{On two theorems of combinatory analysis and some allied identities},
Proc. London Math. Soc. (2) \textbf{16} (1917), 315--336.

\bibitem{SW98}
A. Schilling and S. O. Warnaar,
\textit{Supernomial coefficients, polynomial identities and $q$-series},
The Ramanujan Journal \textbf{2} (1998), 459--494.

\bibitem{SW00}
A. Schilling and S. O. Warnaar,
\textit{A generalization of the $q$-Saalsch\"utz sum and the Burge transform},
in \textit{Physical Combinatorics}, pp. 163--183, 
M.~Kashiwara and T.~Miwa eds.,
Progr. Math. \textbf{191}, (Birkh\"auser, Boston, 2000).
 
\end{thebibliography}

\end{document}